\documentclass[12pt]{article}
\usepackage{e-jc}
\usepackage{latexsym, amsmath, amsthm, amssymb,amscd}
\usepackage{multirow}
\usepackage{array}
\usepackage[dvips]{graphicx}
\usepackage{pict2e}
\usepackage{pstricks}
\newtheorem{thm}{Theorem}[section]
\newtheorem{rmk}{Remark}[section]

\newtheorem{defn}[thm]{Definition}
\newtheorem{lem}[thm]{Lemma}

\numberwithin{equation}{section}
\def\pf{\noindent {\it Proof.} }

\def\stat{{\mathrm stat}}
\def\st{{\mathrm st}}
\def\inv{{\mathrm inv}}
\def\los{{\mathrm los}}
\def\cro{{\mathrm cr}}
\def\crol{{\mathrm cr^{(\ell)}}}
\def\croc{{\mathrm cr^{(c)}}}
\def\ov{{\mathrm ovl}}
\def\emb{{\mathrm emb}}
\def\occ{{\mathrm occ}}
\def\occs{{\mathrm \occ_{\sigma}}}

\def\1{{\bf 1}}
\def\2{{\bf 2}}
\def\T{{\mathcal T}}

\def\mafig{\psline(1;0)(0.8;5)(1;10)(0.8;15)(1;20)(0.8;25)(1;30)(0.8;35)(1;40)(0.8;45)(1;50)(0.8;55)(1;60)}

\def\casA{ \pscircle(0,0){1}
 \uput[ur](0,1){\scriptsize 1}
 \psdots(0,1)(1,0)(0,-1)(-1,0)
 \psline[linewidth=.04](0,1)(0,-1)
 \psline[linewidth=.04](1,0)(-1,0)}

\def\casB{ \pscircle(0,0){1}
 \uput[ur](0,1){\scriptsize 1}
 \uput[dl](0,-1){\scriptsize $i$}
 \psdots(0,1)(0,-1)
 \psline[linewidth=.04](0,1)(0,-1)
 \SpecialCoor
 \psline[linewidth=.04](1;30)(1;150)
 \psline[linewidth=.04](1;-30)(1;210)
 \rput{-30}(0,0){\mafig}
 \rput{150}(0,0){\mafig}}

 \def\casCa{
 \pscircle(0,0){1}
 \uput[ur](0,1){\scriptsize 1}
 \uput[r](0.86,-0.5){\scriptsize $i$}
 \uput[l](-0.86,0.5){\scriptsize $n$}
 \SpecialCoor
 \psdots(0,1)(1;-30)(1;150)
 \psline[linewidth=.04](1;-30)(1;150)
 \psline[linewidth=.04](1;90)(1;-150)
 \psline[linewidth=.04](1;30)(1;-90)
 \rput{30}(0,0){\mafig}
 \rput{-150}(0,0){\mafig}}

 \def\casCb{\pscircle(0,0){1}
 \uput[ur](0,1){\scriptsize 1}
 \uput[r](1,0){\scriptsize $i$}
 \uput[l](-1,0){\scriptsize $j$}
 \SpecialCoor
 \psdots(0,1)(1;0)(1;180)
 \psline[linewidth=.04](1;0)(1;180)
 \psline[linewidth=.04](1;60)(1;-60)
 \psline[linewidth=.04](1;120)(1;-120)
 \rput{60}(0,0){\mafig}
 \rput{-120}(0,0){\mafig}}

 \def\casDa{ \pscircle(0,0){1}
 \uput[ur](0,1){\scriptsize 1}
 \uput[r](1,0){\scriptsize $i$}
 \uput[dr](0.5,-0.86){\scriptsize $j$}
 \uput[l](-1,0){\scriptsize $p$}
 \uput[l](-0.5,0.86){\scriptsize $n$}
 \SpecialCoor
 \psdots(0,1)(1;0)(1;-60)(1;120)(1;180)
 \psline[linewidth=.04](1;90)(1;-150)
 \psline[linewidth=.04](1;0)(1;120)
 \psline[linewidth=.04](1;-60)(1;180)
 \psline[linewidth=.04](1;30)(1;-90)
 \rput{30}(0,0){\mafig}
 \rput{-150}(0,0){\mafig}
 \rput{-60}(0,0){\mafig}
 \rput{120}(0,0){\mafig}}

 \def\casDb{\pscircle(0,0){1}
 \uput[ur](0,1){\scriptsize 1}
 \uput[r](0.86,0.5){\scriptsize $i$}
 \uput[r](0.86,-0.5){\scriptsize $j$}
 \uput[l](-0.86,-0.5){\scriptsize $p$}
 \uput[l](-0.86,0.5){\scriptsize $q$}
 \SpecialCoor
 \psdots(0,1)(1;30)(1;-30)(1;150)(1;-150)
 \psline[linewidth=.04](1;60)(1;-60)
 \psline[linewidth=.04](1;120)(1;-120)
 \psline[linewidth=.04](1;30)(1;150)
 \psline[linewidth=.04](1;-30)(1;-150)
 \rput{60}(0,0){\mafig}
 \rput{-120}(0,0){\mafig}
 \rput{-30}(0,0){\mafig}
 \rput{150}(0,0){\mafig}}

\title{Average values of some Z-parameters\\
 in a random set partition}
\author{Anisse Kasraoui\thanks{The author was supported by grant no.\ 090038012 from the Icelandic
Research Fund and  grant S9607-N13 from Austrian Science Foundation FWF in the framework of the National Research Network
``Analytic Combinatorics and Probabilistic Number Theory{}''.}\\
 \small Fakult\"at f\"ur Mathematik, Universit\"at Wien\\[-0.8ex]
\small Nordbergstra\ss e 15, A-1090 Vienna, Austria\\
\small\texttt{anisse.kasraoui@univie.ac.at}}

\date{\dateline{April 6, 2011}{November 21, 2011}\\
   \small Mathematics Subject Classification: 05A18; 05A15, 05A16}

\begin{document}

\maketitle

\begin{abstract}
 We find exact and asymptotic formulas for the average values of several statistics
on set partitions: of Carlitz's $q$-Stirling distributions, of the numbers of crossings 
in linear and circular representations of set partitions, of the numbers of overlappings 
and embracings, and of the numbers of occurrences of a 2-pattern.   
\end{abstract}


\section{Introduction}

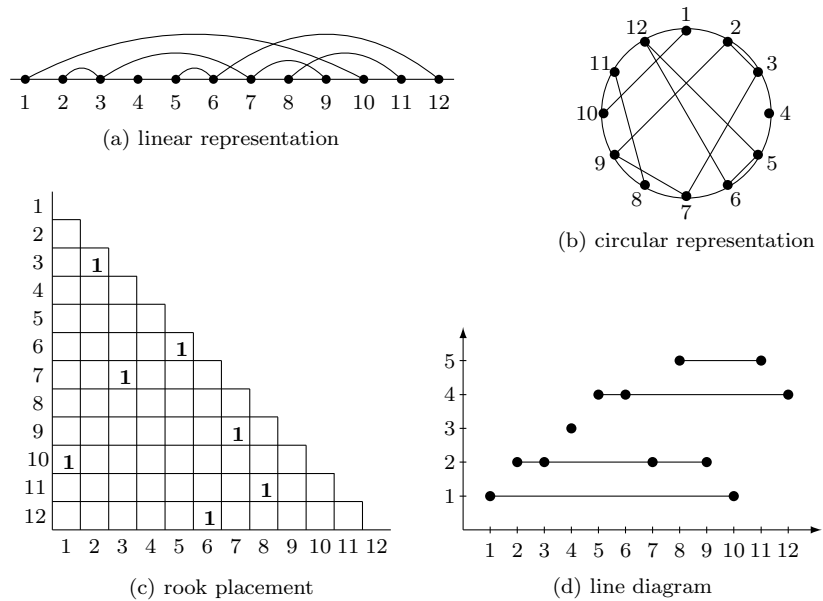
\begin{figure}[t!]\label{fig:representations}
\begin{center}
{\setlength{\unitlength}{1mm}
\begin{picture}(55,10)(0,-10)
\put(-2,0){\line(1,0){59}}
\put(0,0){\circle*{1,2}}\put(0,0){\makebox(0,-6)[c]{\scriptsize 1}}
\put(5,0){\circle*{1,2}}\put(5,0){\makebox(0,-6)[c]{\scriptsize 2}}
\put(10,0){\circle*{1,2}}\put(10,0){\makebox(0,-6)[c]{\scriptsize
3}}
\put(15,0){\circle*{1,2}}\put(15,0){\makebox(0,-6)[c]{\scriptsize
4}}
\put(20,0){\circle*{1,2}}\put(20,0){\makebox(0,-6)[c]{\scriptsize
5}}
\put(25,0){\circle*{1,2}}\put(25,0){\makebox(0,-6)[c]{\scriptsize
6}}
\put(30,0){\circle*{1,2}}\put(30,0){\makebox(0,-6)[c]{\scriptsize
7}}
\put(35,0){\circle*{1,2}}\put(35,0){\makebox(0,-6)[c]{\scriptsize
8}}
\put(40,0){\circle*{1,2}}\put(40,0){\makebox(0,-6)[c]{\scriptsize
9}}
\put(45,0){\circle*{1,2}}\put(45,0){\makebox(0,-6)[c]{\scriptsize
10}}
\put(50,0){\circle*{1,2}}\put(50,0){\makebox(0,-6)[c]{\scriptsize
11}}
\put(55,0){\circle*{1,2}}\put(55,0){\makebox(0,-6)[c]{\scriptsize
12}}
\qbezier(0, 0)(22.5, 12)(45, 0) \qbezier(5, 0)(7.5, 3)(10, 0)
\qbezier(10,0)(20,7)(30,0) \qbezier(30,0)(35,5)(40,0)
\qbezier(20,0)(22.5,3)(25,0) \qbezier(25,0)(40,12)(55,0)
\qbezier(35,0)(42.5,7)(50,0)
\put(26,-5){\makebox(0,-6)[c]{\scriptsize (a) linear
representation}}
\end{picture}
} \hspace{1.5cm}
{ \setlength{\unitlength}{1.1mm}
\begin{picture}(25,25)(-10,-5)
 \put(0, 0){\circle{20,5}} \put(0, 10){\circle*{1,2}} \put(5, 8.66){\circle*{1,2}}
 \put(8.66, 5){\circle*{1,2}} \put(10, 0){\circle*{1,2}} \put(8.66, -5){\circle*{1,2}}
 \put(5, -8.66){\circle*{1,2}} \put(0,-10){\circle*{1,2}} \put(-5, -8.66){\circle*{1,2}}
 \put(-8.66, -5){\circle*{1,2}} \put(-10, 0){\circle*{1,2}} \put(-8.66, 5){\circle*{1,2}}
 \put(-5, 8.66){\circle*{1,2}}
 \put(0, 12){\makebox(0,0)[c]{\scriptsize 1}} \put(6, 10.39){\makebox(0,0)[c]{\scriptsize 2}}
 \put(10.39, 6){\makebox(0,0)[c]{\scriptsize 3}} \put(12,0){\makebox(0,0)[c]{\scriptsize 4}}
 \put(10.39,-6){\makebox(0,0)[c]{\scriptsize 5}} \put(6,-10.39){\makebox(0,0)[c]{\scriptsize 6}}
 \put(0,-12){\makebox(0,0)[c]{\scriptsize 7}} \put(-6,-10.39){\makebox(0,0)[c]{\scriptsize 8}}
 \put(-10.39,-6){\makebox(0,0)[c]{\scriptsize 9}} \put(-12,0){\makebox(0,0)[c]{\scriptsize 10}}
 \put(-10.39,6){\makebox(0,0)[c]{\scriptsize 11}} \put(-6,10.39){\makebox(0,0)[c]{\scriptsize 12}}
 \qbezier(-10,0)(-5, 5)(0,10)
 \qbezier(5,8.66)(6.83, 6.83)(8.66,5) \qbezier(8.66,5)(4.33, -2.5)(0,-10)
 \qbezier(0,-10)(-4.33,-7.5)(-8.66,-5) \qbezier(-8.66,-5)(-1.83,1.83)(5,8.66)
 \qbezier(8.66,-5)(6.83,-6.83)(5,-8.66) \qbezier(5,-8.66)(0,0)(-5,8.66)\qbezier(-5,8.66)(1.83,1.83)(8.66,-5)
 \qbezier(-5, -8.66)(-6.83, -1.83)(-8.66, 5)
\put(0,-15){\makebox(0,-1)[c]{\scriptsize (b) circular
representation}}
\end{picture}}
\end{center}
\begin{center}
{\setlength{\unitlength}{0.75mm}
\begin{picture}(62,60)(-8,0)
\put(0,0){\line(1,0){60}}\put(0,5){\line(1,0){55}}\put(0,10){\line(1,0){50}}
\put(0,15){\line(1,0){45}}\put(0,20){\line(1,0){40}}\put(0,25){\line(1,0){35}}
\put(0,30){\line(1,0){30}}\put(0,35){\line(1,0){25}}\put(0,40){\line(1,0){20}}
\put(0,45){\line(1,0){15}}\put(0,50){\line(1,0){10}}\put(0,55){\line(1,0){5}}
\put(0,0){\line(0,1){60}}\put(5,0){\line(0,1){55}}\put(10,0){\line(0,1){50}}
\put(15,0){\line(0,1){45}}\put(20,0){\line(0,1){40}}\put(25,0){\line(0,1){35}}
\put(30,0){\line(0,1){30}}\put(35,0){\line(0,1){25}}\put(40,0){\line(0,1){20}}
\put(45,0){\line(0,1){15}}\put(50,0){\line(0,1){10}}\put(55,0){\line(0,1){5}}
\put(0,10){\makebox(6,4)[c]{\scriptsize \1}}
\put(5,45){\makebox(6,4)[c]{\scriptsize \1}}
\put(10,25){\makebox(6,4)[c]{\scriptsize \1}}
\put(20,30){\makebox(6,4)[c]{\scriptsize \1}}
\put(25,0){\makebox(6,4)[c]{\scriptsize \1}}
\put(30,15){\makebox(6,4)[c]{\scriptsize \1}}
\put(35,5){\makebox(6,4)[c]{\scriptsize \1}}
\put(0,0){\makebox(5,-6)[c]{\scriptsize $1$}}
\put(5,0){\makebox(5,-6)[c]{\scriptsize $2$}}
\put(10,0){\makebox(5,-6)[c]{\scriptsize $3$}}
\put(15,0){\makebox(5,-6)[c]{\scriptsize $4$}}
\put(20,0){\makebox(5,-6)[c]{\scriptsize $5$}}
\put(25,0){\makebox(5,-6)[c]{\scriptsize $6$}}
\put(30,0){\makebox(5,-6)[c]{\scriptsize $7$}}
\put(35,0){\makebox(5,-6)[c]{\scriptsize $8$}}
\put(40,0){\makebox(5,-6)[c]{\scriptsize $9$}}
\put(45,0){\makebox(5,-6)[c]{\scriptsize $10$}}
\put(50,0){\makebox(5,-6)[c]{\scriptsize $11$}}
\put(55,0){\makebox(5,-6)[c]{\scriptsize $12$}}
\put(0,0){\makebox(-6,5)[c]{\scriptsize $12$}}
\put(0,5){\makebox(-6,5)[c]{\scriptsize $11$}}
\put(0,10){\makebox(-5,5)[c]{\scriptsize $10$}}
\put(0,15){\makebox(-5,5)[c]{\scriptsize $9$}}
\put(0,20){\makebox(-5,5)[c]{\scriptsize $8$}}
\put(0,25){\makebox(-5,5)[c]{\scriptsize $7$}}
\put(0,30){\makebox(-5,5)[c]{\scriptsize $6$}}
\put(0,35){\makebox(-5,5)[c]{\scriptsize $5$}}
\put(0,40){\makebox(-5,5)[c]{\scriptsize $4$}}
\put(0,45){\makebox(-5,5)[c]{\scriptsize $3$}}
\put(0,50){\makebox(-5,5)[c]{\scriptsize $2$}}
\put(0,55){\makebox(-5,5)[c]{\scriptsize $1$}}
\put(30,-9){\makebox(0,-3)[c]{\scriptsize (c) rook placement}}
\end{picture}}
\hspace{1cm}
{\setlength{\unitlength}{0.9mm}
\begin{picture}(54,50)(0,0)
\put(0,0){\vector(1,0){53}}
\put(4,0){\makebox(0,-6)[c]{\scriptsize$1$}}\put(8,0){\makebox(0,-6)[c]{\scriptsize$2$}}
\put(12,0){\makebox(0,-6)[c]{\scriptsize$3$}}\put(16,0){\makebox(0,-6)[c]{\scriptsize$4$}}
\put(20,0){\makebox(0,-6)[c]{\scriptsize$5$}}\put(24,0){\makebox(0,-6)[c]{\scriptsize$6$}}
\put(28,0){\makebox(0,-6)[c]{\scriptsize$7$}}\put(32,0){\makebox(0,-6)[c]{\scriptsize$8$}}
\put(36,0){\makebox(0,-6)[c]{\scriptsize$9$}}\put(40,0){\makebox(0,-6)[c]{\scriptsize$10$}}
\put(44,0){\makebox(0,-6)[c]{\scriptsize$11$}}\put(48,0){\makebox(0,-6)[c]{\scriptsize$12$}}
\put(4,-0.5){\line(0,1){1}}\put(8,-0.5){\line(0,1){1}}
\put(12,-0.5){\line(0,1){1}}\put(16,-0.5){\line(0,1){1}}
\put(20,-0.5){\line(0,1){1}}\put(24,-0.5){\line(0,1){1}}
\put(28,-0.5){\line(0,1){1}}\put(32,-0.5){\line(0,1){1}}
\put(36,-0.5){\line(0,1){1}}\put(40,-0.5){\line(0,1){1}}
\put(44,-0.5){\line(0,1){1}}\put(48,-0.5){\line(0,1){1}}
\put(0,0){\vector(0,1){30}}
\put(0,5){\makebox(-4,0)[c]{\scriptsize$1$}}\put(0,10){\makebox(-4,0)[c]{\scriptsize$2$}}
\put(0,15){\makebox(-4,0)[c]{\scriptsize$3$}}\put(0,20){\makebox(-4,0)[c]{\scriptsize$4$}}
\put(0,25){\makebox(-4,0)[c]{\scriptsize$5$}}
\put(-0.5,5){\line(1,0){1}}\put(-0.5,10){\line(1,0){1}}
\put(-0.5,15){\line(1,0){1}}\put(-0.5,20){\line(1,0){1}}
\put(-0.5,25){\line(1,0){1}}
\put(4,5){\circle*{1.5}}\put(40,5){\circle*{1.5}}
\put(8,10){\circle*{1.5}}\put(12,10){\circle*{1.5}}\put(28,10){\circle*{1.5}}\put(36,10){\circle*{1.5}}
\put(16,15){\circle*{1.5}}
\put(20,20){\circle*{1.5}}\put(24,20){\circle*{1.5}}\put(48,20){\circle*{1.5}}
\put(32,25){\circle*{1.5}}\put(44,25){\circle*{1.5}}
\put(4,5){\line(1,0){36}}\put(8,10){\line(1,0){28}}\put(20,20){\line(1,0){28}}\put(32,25){\line(1,0){12}}
\put(25,-7){\makebox(0,-3)[c]{\scriptsize (d) line diagram}}
\end{picture}
}
\end{center}
\vspace{1.0cm} \caption{{\small Some
representations of the set partition $\pi=1\,10/2\,3\,7\,9/4/5\,6\,12/8\,11$}}
\end{figure}

  Among all basic combinatorial structures, set partitions are probably one
of the most  attractive. Like their close cousin, the
permutations, they have a remarkably rich combinatorial structure.
This is partly due to the fact that they can be represented in many equivalent
way: in terms of words, functions, graphs, placements of non-attacking
rooks in a triangular diagram, line diagrams, etc. (see
Figure~\ref{fig:representations} for some examples).  
During the past decades several
 statistics (here, statistic means an integer valued combinatorial parameter), 
also called parameters, on these different representations of set partitions have been
introduced and studied. For instance, natural statistics on set
partitions are the number of inversions in restricted growth
functions or the numbers of crossings of arcs in linear and circular
representations of set partitions. In very rare cases, one can
determine the exact distribution of the statistic (i.e., for each
$n$, the number of set partitions of $\{1,2,\ldots,n\}$ which
have a certain value of the considered statistic), but most often,
even though one has good information on the parameter, such as,
for example, in form of an explicit formula for the
generating function, there is no exact expression for
the distribution. 
 The by far most important  descriptive quantity associated to a statistic is 
its \emph{average value} (also called \emph{mean value} or
\emph{expectation}).  
The purpose of this paper is to compute the 
average values of several statistics on set partitions.
Some of these statistics have been introduced in the literature previously and others are new. 

\section{Background and main results}

\subsection{Basic facts}

The central objects of this paper are set partitions.  A partition
$\pi$ of a set $S$ is a collection of nonempty and
mutually disjoint sets $\pi=\{A_1,A_2,\ldots,A_k\}$ whose union
is~$S$. The sets $A_i$ are called the \emph{blocks} of the
partition, and a partition into~$k$ blocks is called
a \emph{$k$-partition}. Let us denote by $\Pi(S)$ the family of
all partitions of $S$ and by $\Pi^k(S)$ the family of all
$k$-partitions of $S$.
 For instance, $\Pi^2(\{a,b,c,d\})$ consists of the elements:
\begin{align*}
&\{\{a\},\{b,c,d\}\}\quad\{\{b\},\{a,c,d\}\}\quad\{\{c\},\{a,b,d\}\}
\quad\{\{d\},\{a,b,c\}\}\\
&\{\{a,b\},\{c,d\}\}\quad\{\{a,c\},\{b,d\}\}\quad\{\{a,d\},\{b,c\}\}
\end{align*}
 The number $S_{n,k}$ of $k$-partitions of an $n$-set is called
\emph{Stirling number of the second kind}, while the number $B_n$
of all partitions of an $n$-set is called \emph{Bell number}.
Hence, by definition, we have $B_n= \sum_{k=1}^n S_{n,k}$.

 Almost all statistics considered on set partitions require that the
set~$S$ to be partitioned  is linearly ordered. Without loss of
generality, from now on, we will consider only partitions of subsets
of the set $\mathbb{P}=\{1,2,\ldots\}$ of positive integers. With
$[n]:=\{1,2,\ldots,n\}$, we will set $\Pi_n:=\Pi([n])$, 
$\Pi_n^k:=\Pi^k([n])$ and $\Pi:=\bigcup_{n\geq1}\Pi_n$.  A convenient way, 
called \emph{standard form} (also called \emph{canonical order}), 
to write a set partition $\pi\in\Pi(S)$, $S\subseteq \mathbb{P}$, 
is to arrange the blocks of $\pi$ in increasing order
of their minima and each block in increasing order. For instance,
the partition $\pi$ of $[9]$ consisting of the five blocks 
$\{1,4,7\},\, \{3,9\},\, \{2\},\, \{5\}\, \textrm{and}\, \{6,8\} $
is written in standard form as $\pi=1\;4\;7 /2/3\;9/5 /6\;8$.\\

 A substantial number of combinatorial statistics
on set partitions were defined in terms of the {\emph restricted
growth function} representation for set partitions. It is well-known
that~$k!S_{n,k}$ is the number of surjections from an $n$-set to a
$k$-set. Actually, Stirling numbers also count a class of functions,
called restricted growth functions. A word (or equivalently, a
mapping)  $w=w_1\,w_2\cdots\,w_n$ is said to  be a restricted growth
function (RGF for short) of length $n$ and size $k$  if:
\begin{enumerate}
\item $w_1=1$,
\item $w_{i+1}\leq \max\{w_{1},w_{2},\ldots, w_{i}\}+1$ for $1\leq i\leq n-1$,
\item $\max\{w_{1},w_{2},\ldots, w_{n}\}=k$.
\end{enumerate}
 The usual way to associate bijectively set partitions of~$[n]$ into~$k$ blocks
 and~RGF of length~$n$ and size~$k$ works as follows: to a partition
$\pi=B_1/B_2/\cdots/B_k$ of~$[n]$ written in standard form, we
associate the RGF $w(\pi)=w_1w_2\ldots w_n$ such that $w_i=j$ if and only
if $i\in B_j$. For instance, under this correspondence,
the partition $\pi=1\;4\;7 /2/3\;9/5 /6\;8$ is sent
to the RGF $w(\pi)=1\,2\,3\,1\,4\,5\,1\,5\,3$.
 The reader who knows the plethora of permutation and word statistics
can imagine the amount of statistics that can be defined on set
partitions via restricted growth functions. For instance, several statistics were
defined in terms of inversions and noninversions, and more
generally, in terms of counting occurrences of patterns in
restricted growth functions.

  Other interesting representations of set partitions suggest 
other combinatorial parameters. We will present some of them
when they are needed.

\subsection{Z-statistics and main results}\label{sect:standmap}

  Computation of average values of combinatorial parameters is often 
entertaining,
but can be a very difficult task in some cases (for instance, consider
the parameter ``length of the longest increasing subsequence'' in permutations).
In this paper, we will consider some statistics whose computation of their (at
least exact) average value is not trivial, but not as difficult as
the one we mentioned previously. For instance, two of the parameters
considered in the present paper are the numbers of crossings in
linear and circular representations of set partitions (see
Figure~\ref{fig:representations}). A little
moment's thought will convince the reader that the computation of
their (at least exact) average value  is not completely trivial.

  Depending on which statistic on set partitions is considered, 
it is more or less difficult to find its average value. 
As an example, let us consider the statistics `number
of blocks' and `number of crossings'. While a simple direct combinatorial
argument gives the average number of blocks in a random set partition of $[n]$ (it
is equal to $\frac{B_{n+1}}{B_{n}}-1$), a moment's thought will
convince the reader that the computation of the average number of
crossings in the circular representation of a random set partition
of $[n]$, which, as we shall see later, is actually equal to
$$
\frac{n}{2}\frac{B_{n+1}}{B_{n}}+\frac{3n}{2}-\frac{n(4n+1)}{2}\frac{B_{n-1}}{B_{n}}
-{n\choose 2}\frac{B_{n-2}}{B_{n}}+{n\choose 4}\frac{B_{n-4}}{B_{n}}
$$
is not completely trivial.

 Perhaps the most natural and powerful approach to deal with average
values is to use the linearity of the expectation by finding ``nice
decompositions'' of the statistics we consider. This was done in the
author's thesis~\cite{Ka} for the average number of crossings in the
linear representation of a random set partition. It turns out that the
decomposition used there is shared by a large number of existing
statistics on set partitions and that their average values can be
obtained almost ``mechanically''. In order to make these remarks more
precise, we need additional terminology.

  Given a subset $S\subseteq \mathbb{P}$ with cardinality $|S|=n$, the \emph{standardization map} $\st_S$ is the (unique)
order-preserving bijection $\st_S: S \to [n]$. When the subset $S$ is clear from the context, we will drop the subscript. We let $\st_S$
act element-wise on objects built using $S$ as label. For instance, the set partition $\pi=2\,9/4\,10/5/7\,11/8$
of~$S=\{2,4,5,7,8,9,10,11\}$ is sent after standardization to the set partition $\st(\pi)=1\,6/2\,7/3/4\,8/5$.

\begin{defn}\label{defn:splittable}
  A set partition statistic $STAT$ is said to be a \emph{Z-statistic}
if for any $\pi\in\Pi$,
\begin{align}
  STAT(\pi)=\sum\limits_{\substack{A,B\in \pi \\
A\neq B}} STAT \left( \st(\{A,B\}) \right).\label{eq:defn_splittable}
\end{align}
\end{defn}

  The terminology is inspired from a paper of Galovich and White~\cite{GaWh},
but originally, the term $Z$-statistic is for a famous statistic introduced 
by Bressoud and Zeilberger~\cite{BrZe}. It is worth noting that via the RGF representation, 
our \emph{Z-statistics} on set partitions could be seen as  $Z$-statistics 
on words (with the definition of Galovitch and White).

For the present paper, the interest of Z-statistics on set partitions 
is that the computation of their average value
in a random partition is essentially equivalent to the computation
of their average value in a random 2-partition (see Theorem~\ref{thm:GfunMean}), 
the latter being often easier to calculate.

An important number of existing statistics on set partitions are in fact Z-statistics.
In this paper, we will compute average values of Carlitz's $q$-Stirling distributions, 
the numbers of crossings in linear and circular representations, overlappings and embracings, 
and the number of occurrences of a 2-pattern in restricted growth functions.


\subsubsection{\bf Carlitz's $q$-Stirling distributions}\label{section:carlitz}

The most common $q$-analogs of Stirling numbers of the second kind
found in the literature are probably the $q$-Stirling numbers
$S_{n,k}(q)$  and $\tilde S_{n,k}(q)$ essentially described by Carlitz and Gould~\cite{Carlitz,Gou}.
They can be defined as follows:
\begin{eqnarray}\label{eq:stirling1}
 S_{n,k}(q)=q^{k-1}\,S_{n-1,k-1}(q)+[k]_q\,S_{n-1,k}(q)\qquad (n\geq k\geq 1),
\end{eqnarray}
and $\tilde S_{n,k}(q)=q^{-{k\choose 2}}\, S_{n,k}(q)$, with
$S_{n,k}(q)=\tilde S_{n,k}(q)=\delta_{n\,k}$ if $n=0$ or $k=0$.

 Our interest in these $q$-Stirling numbers  comes from the fact
that they arise as generating functions for the distribution of
several statistics on set partitions. Such statistics have been
given by Milne~\cite{Mi}, Garsia and Remmel~\cite{GaRe}, Sagan~\cite{Sag},
Wachs and White~\cite{WW}, White~\cite{Wh}, Leroux~\cite{Ler}, de M\'edicis and Leroux~\cite{MedLer},
Ehrenborg and Readdy~\cite{EhRe}, Ksavrelof and Zeng~\cite{KsZe}. 
This list is probably not exhaustive. It is natural to ask for the 
expectation of these ``$q$-Stirling distributions''.

We will say that a parameter  $\stat$  defined on set partitions has
the $q$-Stirling distribution~$S(q)$ ($\tilde S(q)$, respectively) if its
generating function on each $\Pi_n^k$ is given by $S_{n,k}(q)$ 
($\tilde S_{n,k}(q)$, respectively), i.e., for any $n\geq k\geq 1$, we have
\begin{align*}
\sum_{\pi\in\Pi_n^k} q^{\stat(\pi)}=S_{n,k}(q) \;  \textrm{($=\tilde S_{n,k}(q)$, respectively)}.
\end{align*}

 \begin{thm}\label{thm:Mean_Qstirling}
Let $\stat$ be a parameter with the $q$-Stirling distribution~$S(q)$.
Then the average value $\mu_n$ of $\stat$  in a random set partition
of $[n]$ satisfies 
\begin{align}
\mu_{n}&=-\frac{1}{4}\frac{B_{n+2}}{B_n}+\left(\frac{n}{2}+\frac{1}{4}\right)\frac{B_{n+1}}{B_n}
+\left(-\frac{n}{2}+\frac{1}{4}\right)\label{eq:ExactMean_qstirling}\\
      &=\frac{n^2}{2\log n}\left( 1+\frac{\log\log n}{\log
n}\left(1+o\left(1\right)\right) \right) \quad\text{as $n\rightarrow\infty$.} \label{eq:AsymptoticMean_qstirling}
\end{align}
 The average value $\mu_{n,k}$ of $\stat$  in a random set
partition of $[n]$ into $k$~blocks satisfies 
\begin{align}
\mu_{n,k}&=\frac{1}{2} n(k - 1)- \frac{1}{2}{k\choose 2}
  + \frac{1}{2}(n+1-k)\frac{S_{n,k-1}}{S_{n,k}}\label{eq:ExactMeanBlock_qstirling}\\
 &=\frac{1}{2} n(k - 1)- \frac{1}{2}{k\choose 2}+ o(1)\quad\text{as $n\rightarrow\infty$.}\label{eq:AsymptoticMeanBlock_qstirling}
\end{align}
\end{thm}

Note that it is almost immediate to derive from the
above results the expectations of the $q$-Stirling
distribution~$\tilde S(q)$, which we
denoted by $\tilde\mu_{n}$ and $\tilde\mu_{n,k}$, respectively. 
Indeed, by the relation $\tilde S_{n,k}(q)=q^{-{k\choose 2}}\, S_{n,k}(q)$, we have
\begin{align*}
\tilde\mu_{n,k}=\mu_{n,k}-{k\choose2}
\quad\text{and}\quad\tilde\mu_{n}=\mu_{n}-\frac{1}{B_n}\sum_{k=1}^n{k\choose2}S_{n,k},
\end{align*}
and thus
\begin{align*}
\tilde\mu_{n,k}&=\frac{1}{2} n(k - 1)- \frac{3}{2}{k\choose 2}
  + \frac{1}{2}(n+1-k)\frac{S_{n,k-1}}{S_{n,k}}\\
\tilde\mu_{n}&=-\frac{3}{4}\frac{B_{n+2}}{B_n}+\left(\frac{n}{2}+\frac{7}{4}\right)\frac{B_{n+1}}{B_n}
-\left(\frac{n}{2}+\frac{1}{4}\right),
\end{align*}
where the last equality follows from the identity 
$\sum_{k=1}^n k(k-1)S_{n,k}=B_{n+2}-3B_{n+1}+B_n$ 
(for a proof, see e.g.~\cite{Sach}).\\[-0.3cm]

Almost all known statistics with the $q$-Stirling distribution $S(q)$ 
are in fact Z-statistics. In order to prove Theorem~\ref{thm:Mean_Qstirling}, we will work 
with the first known parameter with the distribution $S(q)$: this is the statistic $\los$
(originally denoted $I^{m}$), first introduced by Milne (see~\cite[Section~4]{Mi}) 
to answer a question of Garsia. The parameter $\los$ can also be defined by using the 
non-attacking rook placement model for set partitions (see e.g.~\cite{WW}).                   

It is not hard to see that $\los$ can be defined for a 
set partition $\pi=B_1/B_2/\ldots/B_k$ written in standard form 
by $$\los(\pi)=|B_2|+2 |B_3|+\cdots+(k-1)|B_k|.$$
 It is easy to see that the parameter $\los$ is a Z-statistic. 
Indeed, for a 2-partition $B_1/B_2$ written in standard form,
we have $\los(B_1/B_2)=|B_2|$ whence for $\pi=B_1/B_2/\ldots/B_k$ we have
\begin{align*}
 \sum_{1\leq i <j\leq k} STAT \left( \st(B_i/B_j) \right)
= \sum_{1\leq i <j\leq k}|B_j|
=\sum_{j=2}^k  (j-1)\, |B_j|=\los(\pi).
\end{align*}


\subsubsection{\bf Number of linear crossings}

One of the most widely occurring property of set partitions is
the ``noncrossing'' property. 
This is  essentially due to the important number of works on noncrossing partitions  which appear
 in different fields of mathematics.  The terminology originated probably 
from the graphical representation 
of set partitions (see e.g.~\cite{Si}). The \emph{linear representation} of
 a partition of a set $S \subseteq \mathbb{P}$ is obtained as follows: draw the elements of
 $S$ in increasing order on a line and join successive  elements of each block by arcs drawn
 in the upper half-plane.

 Alternatively, the elements may be represented on a circle and circularly successive elements of
 each block are joined by chords: this is the \emph{circular representation}. 
For instance, the linear and circular representations  of the partition~$\pi=1\,10/2\,3\,7\,9/4/5\,6\,12/8\,11$  are given
 in Figure~\ref{fig:representations}. Linear and circular representations of set partitions suggest many statistics.
 The most studied is the number of crossings, but other statistics like nestings, alignments,
 $k$-distant crossings, $k$-crossings, $k$-nestings, etc.\ have also appeared in the literature (see e.g.~\cite{ChenSt,DraKim,KaZe}).

   For a set partition $\pi\in\Pi$, we will denote by $\crol(\pi)$ the number of crossings in the linear representation of
 $\pi$. We will refer to  such crossings as  \emph{linear crossings},
in order to distinguish these from the circular case.
Formally, a linear crossing is a sequence $(i_1,j_1),(i_2,j_2)$ of edges of $\pi$,
 such that $i_1<i_2<j_1<j_2$. Here, an edge of $\pi$ is just a pair $(i,j)$ of \emph{consecutive}
elements $i<j$ in the same block $B$, where consecutive means that
there is no element $p\in B$ such that $i<p<j$. For instance, if $\pi$ is the 
set partition represented in Figure~\ref{fig:representations}, then $\crol(\pi)=4$. 
   In order to show that $\crol$ is a Z-statistic, it suffices essentially to 
observe that
a linear crossing is composed of two edges from distinct blocks (details are left to the reader). 

 \begin{thm}\label{thm:Mean_CrossingLinear}
The average number $\mu_n$ of linear crossings in a random set
partition of~$[n]$ satisfies
\begin{align}
\mu_n&=-\frac{5}{4}\frac{B_{n+2}}{B_{n}}+\left(\frac{n}{2}+\frac{9}{4}\right)\frac{B_{n+1}}{B_{n}}
+\frac{n}{2}+\frac{1}{4}\label{eq:ExactMean_crol}\\
      &=\frac{n^2}{2\log n}\left( 1+\frac{\log\log n}{\log
n}\left(1+o\left(1\right)\right) \right)\quad\text{as $n\rightarrow\infty$.}\label{eq:AsymptoticMean_crol}
\end{align}
The average number $\mu_{n,k}$ of linear crossings in a random set
partition of~$[n]$ into~$k$ blocks satisfies
\begin{align}
\mu_{n,k}&=\frac{1}{2}n(k-1)-\frac{5}{2}{k\choose 2}
+\frac{3}{2}(n+1-k)\frac{S_{n,k-1}}{S_{n,k}}\label{eq:ExactMeanBlock_crol}\\
         &=\frac{1}{2}n(k-1)-\frac{5}{2}{k\choose 2}+ o(1)\quad\text{as $n\rightarrow\infty$.}
\label{eq:AsymptoticMeanBlock_crol}
\end{align}
\end{thm}


\subsubsection{\bf Number of circular crossings}

 For a set partition $\pi\in\Pi$, we will denote by $\croc(\pi)$ the number of crossings in the circular representation of
 $\pi$. We will refer to  such crossings as  \emph{circular crossings}.
Formally, a circular crossing is a sequence $(i_1,j_1),(i_2,j_2)$ of
``circular edges'' of $\pi$, such that ${i_1<i_2<j_1<j_2}$. Here, a
circular edge of $\pi$ is a pair $(i,j)$ such that $i,j$ are
consecutive elements in the same block of $\pi$ or $i=\min(B)$ and
$j=\max(B)$ for some block $B\in\pi$.  For instance, if $\pi$ is the set partition 
represented in Figure~\ref{fig:representations}, 
then $\croc(\pi)=9$. We leave the verification 
that $\croc$ is a Z-statistic to the reader.

\begin{thm}\label{thm:Mean_CrossingCircular}
The average number $\mu_n$ of circular crossings in a random set
partition of~$[n]$ satisfies
\begin{align}
\mu_n&=\frac{n}{2}\frac{B_{n+1}}{B_{n}}+\frac{3n}{2}-\frac{n(4n+1)}{2}\frac{B_{n-1}}{B_{n}}
-{n\choose 2}\frac{B_{n-2}}{B_{n}}+{n\choose 4}\frac{B_{n-4}}{B_{n}}\label{eq:ExactMean_croc}\\
    &=\frac{n^2}{2\log n}\left( 1+\frac{\log\log n}{\log
n}\left(1+o\left(1\right)\right) \right)\quad\text{as $n\rightarrow\infty$.}\label{eq:AsymptoticMean_croc}
\end{align}
The average number $\mu_{n,k}$ of circular crossings in a random set
partition of $[n]$ into $k$ blocks satisfies
\begin{align}
\mu_{n,k}&=\frac{1}{2}n(k-1)-\frac{1}{2}n(4n-5k+1)\frac{S_{n-1,k-1}}{S_{n,k}}
-10{n\choose 2}\frac{S_{n-2,k-2}}{S_{n,k}}\label{eq:ExactMeanBlock_croc}\\
&\;+{n\choose 4}\frac{S_{n-4,k-2}}{S_{n,k}}\nonumber\\
         &=\frac{1}{2}n(k-1)+ o(1)\quad\text{as $n\rightarrow\infty$.}\label{eq:AsymptoticMeanBlock_croc}
\end{align}
\end{thm}


\subsubsection{\bf Number of overlappings}

Flajolet and Schott~\cite{FlS} considered a  special class of set partitions, called \emph{non-overlapping partitions}.
 With the implied order structure, two sets~$B$ and~$B'$  are said to overlap
 if  $\min(B)<\min(B')< \max(B)< \max(B')$. Then, a partition is said
 to be \emph{non-overlapping}  if any two blocks  do not overlap. This property of non-overlapping
could be seen graphically by considering the  line diagram (see
Figure~\ref{fig:representations}) of a set partition. This suggests to
consider statistics like the numbers of overlappings, embracings (``nestings of blocks''),
etc. 

 The number of overlappings of a set partition~$\pi$ will be
denoted~$\ov(\pi)$. For instance, if $\pi$ is the set partition 
represented in Figure~\ref{fig:representations}, 
then $\ov(\pi)=4$.  We leave the verification that~$\ov$ is
a Z-statistic to the reader.

\begin{thm}\label{thm:exactmean_ovl}
The average number $\mu_n$ of overlappings in a random set partition
of $[n]$ satisfies
\begin{align}
\mu_n&=\frac{1}{4}\frac{B_{n+2}}{B_{n}}+\frac{3}{4}\frac{B_{n+1}}{B_{n}}
-\left(n+\frac{5}{4}\right)-\frac{n}{2}\frac{B_{n-1}}{B_{n}}\label{eq:ExactMean_ov}\\
    &=\frac{1}{4}\left(\frac{n}{\log n}\right)^2 \left(
1+2\frac{\log\log n}{\log n}\left(1+o\left(1\right)\right) \right)\quad\text{as $n\rightarrow\infty$}.\label{eq:AsymptoticMean_ov}
\end{align}

The average number $\mu_{n,k}$ of overlappings  in a random set
partition of $[n]$  into $k$ blocks satisfies
\begin{align}
\mu_{n,k}
&=\frac{1}{2}{k\choose2}+n(k-1)\frac{S_{n-1,k-1}}{S_{n,k}}
-\frac{3}{2}(n+k-1)\frac{S_{n,k-1}}{S_{n,k}}\label{eq:ExactMeanBlock_ov}\\
&=\frac{1}{2}{k\choose2}+ o(1)\quad\text{as $n\rightarrow\infty$}.\label{eq:AsymptoticMeanBlock_ov}
\end{align}
\end{thm}


\subsubsection{\bf Number of embracings}
In this section, we consider a natural partner for overlappings. Say
that two sets~$B$ and~$B'$ \emph{embrace} if, with the implied order
structure, $\min(B)<\min(B')\leq \max(B')< \max(B)$. The number of
embracings of a set partition~$\pi$ will be denoted~$\emb(\pi)$ and
 can be easily computed via the line diagram representation. For instance, if $\pi$ is the set partition 
represented in Figure~\ref{fig:representations}, 
then $\emb(\pi)=4$.  We
leave the verification that~$\emb$ is a \linebreak Z-statistic to the reader.

\begin{thm}\label{thm:exactmean_emb}
The average number $\mu_{n}$ of embracings  in a random set
partition of $[n]$ satisfies
\begin{align}
\mu_{n}&=\frac{1}{4} \frac{B_{n+2}}{B_n}-\frac{5}{4}
\frac{B_{n+1}}{B_n}+\frac{3}{4}+\frac{n}{2}\frac{B_{n-1}}{B_n}\label{eq:ExactMean_emb}\\
&=\frac{1}{4}\left(\frac{n}{\log n}\right)^2 \left(
1+2\frac{\log\log n}{\log n}\left(1+o\left(1\right)\right) \right)\quad\text{as $n\rightarrow\infty$}.\label{eq:AsymptoticMean_emb}
\end{align}
The average number $\mu_{n,k}$ of embracings  in a random set
partition of $[n]$  into $k$ blocks satisfies
\begin{align}
\mu_{n,k} &=\frac{1}{2}{k\choose2}-\frac{1}{2}
(k-1)\frac{S_{n,k-1}}{S_{n,k}}+\frac{n}{2}\frac{S_{n-1,k-2}}{S_{n,k}}\label{eq:ExactMeanBlock_emb}\\
&=\frac{1}{2}{k\choose2}+ o(1)\quad\text{as $n\rightarrow\infty$}.\label{eq:AsymptoticMeanBlock_emb}
\end{align}
\end{thm}


\subsubsection{\bf Number of occurrences of a 2-pattern}

 We say that a word $\sigma=\sigma_1\sigma_2\ldots\sigma_r$  in $\{1,2\}^*$
is a $2$-pattern if it contains the two letters~$1$ and~$2$. An
occurrence of the 2-pattern $\sigma$ in a word $w=w_1w_2\ldots
w_n\in\mathbb{P}^{*}$ is an $r$-tuple $(i_1,i_2,\ldots,i_r)$ such
that $\st(w_{i_1}w_{i_2}\ldots w_{i_r})=\sigma$, where $\st$ is the
standardization map (see Section~\ref{sect:standmap}). For instance, 
in the word $w=2\,4\,3\,1\,2\,2\,6\,7\,3\,7$, the $3$-tuples $(1,3,9)$,
$(4,5,6)$ and $(7,8,10)$ are  occurrences of the pattern $1\,2\,2$,
while the $3$-tuples $(2,3,9)$ and $(3,5,6)$ are  occurrences of the
pattern $2\,1\,1$. The number of occurrences of a pattern $\sigma$
in a word $w$ will be denoted $\occs(w)$. Since a set partition can
be identified with its restricted growth function, we will set
$\occs(\pi)=\occs(w(\pi))$, where $w(\pi)$ is, as usual, the RGF
associated to $\pi$.  We leave the verification
that $\occ_{\sigma}$ is a Z-statistic for any $2$-pattern~$\sigma$
to the reader.

\begin{thm}\label{thm:exactmean_occ}
Let $\sigma=\sigma_1\sigma_2\ldots\sigma_r$ be a $2$-pattern of length $r$, and let
$\mu_{n}^{(\sigma)}$ denote the average number of occurrences of the
pattern $\sigma$ in a random set partition of $[n]$.\\

 (i) For fixed $n$, the average value $\mu_{n}^{(\sigma)}$ depends only on the first letter of  $\sigma$.\\

 For $i\in\{1,2\}$, denote by $\mu_{n}^{(\sigma_1=i)}$ the common value of $\mu_{n}^{(\tau)}$
for patterns $\tau$ (of length $r$) satisfying $\tau_1=i$.

(ii) For any integer $n\geq 1$, we have
\begin{align}
 \mu_{n}^{(\sigma_1=1)}+\mu_{n}^{(\sigma_1=2)}
&={n \choose r}
\frac{B_{n+2-r}-B_{n+1-r}}{B_{n}}.\label{eq:ConnectingRelation_occs}
\end{align}

 For $j=0,\ldots,r$, set
\begin{align}
p_{j}(n)&= \frac{(-1)^{r-j}}{2^{r-j+1}}\left({n\choose
j}+\frac{1}{2}{n\choose j-1}\right).\label{eq:defP_j}
\end{align}

(iii) For any integer $n\geq 0$, we have
\begin{align*}
\mu_{n}^{(\sigma_1=2)}
&=\sum_{j=0}^{r}p_{j}(n)\;\frac{B_{n+2-j}}{B_n}\;
-\frac{1}{2}{n\choose
r}\frac{B_{n+1-r}}{B_n}+\frac{(-1)^{r+1}}{2^{r+1}}\\
\mu_{n}^{(\sigma_1=1)}
&=-\sum_{j=0}^{r-1}p_{j}(n)\;\frac{B_{n+2-j}}{B_n}\;
+\frac{1}{2}\left({n\choose r}-\frac{1}{2}{n\choose r-1}\right)
\frac{B_{n+2-r}}{B_n}\\
&\qquad-\frac{1}{2}{n\choose
r}\frac{B_{n+1-r}}{B_n}+\frac{(-1)^{r}}{2^{r+1}}.
\end{align*}

(iv) For any  2-pattern $\sigma$ of length $r$, we have the
asymptotic approximation
\begin{align}\label{eq:AsymptoticMean_occ}
\mu_n^{(\sigma)}&=\frac{1}{2\,r!} n^{2} \left(\log n\right)^{r-2}
\left( 1-(r-2) \frac{\log\log n}{\log n} +o\left(\frac{\log\log
n}{\log n}\right) \right),\quad n\rightarrow\infty.
\end{align}

(v) For the 2-patterns $\sigma$ of length 2, i.e., $\sigma\in\{1\,2\,,2\,1\}$,
we have the more precise asymptotic approximation
\begin{align}\label{eq:AsymptoticMean_occ_length2}
\mu_{n}^{(\sigma)}&=\frac{n^2}{4}\left( 1-\frac{1}{\log
n}-\frac{\log\log n}{\left(\log
n\right)^2}\left(1+o\left(1\right)\right) \right),\quad
n\rightarrow\infty.
\end{align}
\end{thm}

\vspace{0.3cm}

 To illustrate our result, we give the values of $\mu_n^{(\sigma)}$ for 2-patterns~$\sigma$
of length $r\leq 3$.
 For patterns of length $r=2$, we have
\begin{align}
\mu_n^{(2\,1)}&=\frac{1}{8}\frac{B_{n+2}}{B_{n}}
-\left(\frac{n}{4}+\frac{1}{8}\right)\frac{B_{n+1}}{B_{n}}
+\left(\frac{n^2}{4}-\frac{1}{8}\right)-\frac{1}{2} {n \choose 2}\frac{B_{n-1}}{B_{n}}\label{eq:mean21}\\
\mu_n^{(1\,2)}&=-\frac{1}{8}\frac{B_{n+2}}{B_{n}}
+\left(\frac{n}{4}+\frac{1}{8}\right)\frac{B_{n+1}}{B_{n}}
+\left(\frac{n^2}{4}-\frac{n}{2}+\frac{1}{8}\right)-\frac{1}{2} {n
\choose 2}\frac{B_{n-1}}{B_{n}}.\label{eq:mean12}
\end{align}

For patterns of length $r=3$, we have
\begin{align*}
\mu_n^{(2\,1\,1)}=\mu_n^{(2\,1\,2)}=\mu_n^{(2\,2\,1)}&=-\frac{1}{16}\frac{B_{n+2}}{B_{n}}
+\frac{1}{8}\left(n+\frac{1}{2}\right)\frac{B_{n+1}}{B_{n}}
-\frac{1}{8}\left(n^2-\frac{1}{2}\right)\\
&\qquad+\frac{1}{12}n\left(n-\frac{1}{2}\right)(n-1)\frac{B_{n-1}}{B_{n}}
-\frac{1}{2} {n \choose 3}\frac{B_{n-2}}{B_{n}}\\
\mu_n^{(1\,1\,2)}=\mu_n^{(1\,2\,1)}=\mu_n^{(1\,2\,2)}
&=\frac{1}{16}\frac{B_{n+2}}{B_{n}}
-\frac{1}{8}\left(n+\frac{1}{2}\right)\frac{B_{n+1}}{B_{n}}
+\frac{1}{8}\left(n^2-\frac{1}{2}\right)\\
&\qquad+\frac{1}{12}n(n-1)\left(n-\frac{7}{2}\right)\frac{B_{n-1}}{B_{n}}
-\frac{1}{2} {n \choose 3}\frac{B_{n-2}}{B_{n}}.
\end{align*}

\begin{thm}\label{thm:exactmeanblock_occ}
Let $\sigma=\sigma_1\sigma_2\ldots\sigma_r$ be a $2$-pattern of length $r$, and let
$\mu_{n,k}^{(\sigma)}$ denote the average number of occurrences of the
pattern $\sigma$ in a random set partition of $[n]$ into $k$ blocks. \\

 (i) For fixed $n$ and $k$,  the average value $\mu_{n,k}^{(\sigma)}$ depends only on the first letter of  $\sigma$.\\

 For $i\in\{1,2\}$, denote by $\mu_{n,k}^{(\sigma_1=i)}$ the common value of $\mu_{n,k}^{(\tau)}$
for patterns $\tau$ (of length $r$) satisfying $\tau_1=i$.

(ii) For any integer $n\geq k\geq 1$, we have
\begin{align}\label{eq:ConnectingRelationBlock_occs}
 \mu_{n,k}^{(\sigma_1=1)}+\mu_{n,k}^{(\sigma_1=2)}
&={n \choose r} \frac{S_{n+2-r,k}-S_{n+1-r,k}}{S_{n,k}}.
\end{align}

(iii) For any integer $n\geq k\geq 1$, we have
\begin{align*}
\mu_{n,k}^{(\sigma_1=2)}&=\sum_{j=0}^{r}p_{j}(n)\;\frac{S_{n+2-j,k}}{S_{n,k}}\;
-\frac{1}{2}{n\choose
r}\frac{S_{n+1-r,k}}{S_{n,k}}+\frac{(-1)^{r+1}}{2^{r+1}}\frac{S_{n,k-2}}{S_{n,k}}\\
\mu_{n,k}^{(\sigma_1=1)}
&=-\sum_{j=0}^{r-1} p_{j}(n)\;\frac{S_{n+2-j,k}}{S_{n,k}}\;
+\frac{1}{2}\left({n\choose r}-\frac{1}{2}{n\choose
r-1}\right) \frac{S_{n+2-r,k}}{S_{n,k}}\\
&\quad -\frac{1}{2}{n\choose
r}\frac{S_{n+1-r,k}}{S_{n,k}}+\frac{(-1)^{r}}{2^{r+1}}\frac{S_{n,k-2}}{S_{n,k}},
\end{align*}
where $p_j(n)$ is given by \eqref{eq:defP_j}.\\

(iv) As $n\rightarrow\infty$, we have the following asymptotic approximations
\begin{align}
\mu_{n,k}^{(\sigma_1=2)} &=\sum_{j=0}^r p_{j}(n) k^{2-j}\;
 -\frac{1}{2}{n\choose r}k^{1-r} +o(1)\label{eq:AsymptoticMeanBlock_occ2}\\
\mu_{n,k}^{(\sigma_1=1)}&=-\sum_{j=0}^r p_{j}(n) k^{2-j}\;
 +{n\choose r}k^{2-r}-\frac{1}{2}{n\choose r}k^{1-r} +o(1).\label{eq:AsymptoticMeanBlock_occ1}
\end{align}

In particular, for any  2-pattern $\sigma$ of length $r$, we have
\begin{align*}
\mu_{n,k}^{(\sigma)}&\sim\frac{n^{r}}{2\, k^{r-1}\,r!}  \left(k-1 \right)
\quad\text{as $n\rightarrow\infty$}.
\end{align*}
\end{thm}

\vspace{0.3cm}

\subsection{Plan of the paper}

 In the next section, we will prove the key result of the 
paper (Theorem~\ref{thm:GfunMean}), which
expresses average values of Z-statistics as coefficients
in certain generating functions.
Moreover, there we introduce the tools which are necessary to use our key result. 
The purpose of the other sections (except the last) is to prove the results presented in Section~2.
Finally, we will conclude the paper
with some remarks, notably  that the methodology could be adapted to
other exponential families and to more general statistics.

%
%

\section{Average value of a Z-statistic in a random set partition}
We said in the previous section that the computation of the average value
of a Z-statistic in a random set partition is essentially equivalent to the computation
of its average value in a random 2-partition. 
 Let us make this more precise.
 Given a Z-statistic $\stat$, let $v^{\stat}_{n,2}:=\sum_{\pi\in\Pi_n^2}\stat(\pi)$, and
let $V_2^{\stat}(x)$ be the exponential generating function of the
sequence $\left(v^{\stat}_{n,2}\right)_{n\ge0}$, i.e.,
\begin{align}
 V_2^{\stat}(x)=\sum_{n\geq0} v^{\stat}_{n,2}\,\frac{x^{n}}{n!}=\sum_{n\geq0} \sum_{\pi\in\Pi_n^2}\stat(\pi)\,\frac{x^{n}}{n!}.
\end{align}

 The following theorem is the key result of the paper.

\begin{thm}\label{thm:GfunMean}
Let $\stat$ be  a Z-statistic. Denote by  $\mu_n$ ($\mu_{n,k}$, respectively)  
the average value of the statistic $\stat$ in a random
set partition of~$[n]$ (respectively of~$[n]$ into $k$ blocks). Then we
have
\begin{align}
\mu_n&=\frac{1}{B_n}\,\left[\frac{x^n}{n!}\right]\,V_2^{\stat}(x)B(x)\label{eq:ThmGfunMean}\\
\mu_{n,k}&=\frac{1}{S_{n,k}}\,\left[\frac{x^n}{n!}\right]\,V_2^{\stat}(x)S_{k-2}(x),\label{eq:ThmGfunMeanBlock}
\end{align}
 where  $B_n$ and  $S_{n,k}$  are Bell and Stirling numbers, and $B(x)$ and the $S_{j}(x)$'s  are exponential
 generating functions of Bell and Stirling numbers defined in~\eqref{eq:Stirling_VerticalGF} and~\eqref{eq:Bell_GF}.
\end{thm}

 After coefficient extraction, Theorem~\ref{thm:GfunMean} can
lead 
to exact and asymptotic formulas for average values of Z-statistics. For this purpose, 
we need some tools that we develop in the next subsection. It is worth noting 
that, here, by an exact formula,  we mean a finite expression which involves
elementary functions, and Bell and Stirling numbers. We then
present a process to get the average and asymptotic value of Z-statistics. 
In a particular case (which includes all the statistics we consider in this paper), 
the algorithm is deterministic. We will illustrate the process in the 
third subsection. Finally, we will prove our key result,
Theorem~\ref{thm:GfunMean}, in the last subsection.

\subsection{Basic tools}

\subsubsection{Basic properties of Stirling and Bell numbers}
 Recall that the number of $k$-partitions of an $n$-set is
the Stirling number of the second kind  $S_{n,k}$, and the number of
all partitions of an $n$-set is the Bell number $B_n$, so that
we have $B_n = \sum_{k=1}^n S_{n,k}$. Classifying the $k$-partitions
of an $n$-set as to whether they do or do not contain a given element
yields the recurrence:
\begin{align}
S_{n,k}=S_{n-1,k-1}+ k S_{n-1,k}\quad(n\geq
k\geq1),\label{eq:Stirling_Recurrence}
\end{align}
with $S_{j,0}=S_{0,j}=\delta_{j,0}$.

 As noted previously, the number of surjections from an $n$-set $N$ to a
$k$-set $K$ is equal to $k!\,S_{n,k}$. Then, given a set~$X$,
$|X|=x$, classifying all the maps $f:N\to X$ according to the
cardinality of their image leads to the fundamental identity 
\begin{align*}
x^n=\sum_{k=0}^n{x \choose k} k! S_{n,k}=\sum_{k=0}^n S_{n,k} (x)_k,
\end{align*}
from which, by binomial inversion, we deduce the summation formula
\begin{align}
S_{n,k}=\frac{1}{k!}\sum_{j=0}^k(-1)^{j}{k\choose j}
(k-j)^n.\label{eq:Stirling_SummationFormula}
\end{align}
Here and in the rest of the paper,  $(x)_k$ stands for the $k$-th lower factorial $x(x-1)\cdots(x-k+1)$. 
Equation~\eqref{eq:Stirling_SummationFormula} leads to the 
generating function identity
\begin{align}
S_k(t):=\sum_{n\geq
0}S_{n,k}\frac{t^n}{n!}=\frac{1}{k!}\left(e^{t}-1\right)^k,\label{eq:Stirling_VerticalGF}
\end{align}
from which we deduce the double generating function
\begin{align}
F(t,u):=\sum_{n\geq k\geq 0}S_{n,k}u^{k}\frac{t^{n}}{n!}=\sum_{k\geq
0}S_k(t)u^{k}=e^{u\left(e^{t}-1\right)}.\label{eq:Stirling_DoubleGF}
\end{align}
Setting $u=1$, we recover the exponential generating function of
Bell numbers
\begin{align}
B(t):=\sum_{n\geq
0}B_n\frac{t^{n}}{n!}=e^{e^{t}-1}.\label{eq:Bell_GF}
\end{align}
There are of course many ways to recover all the previous identities.
For instance, one could first obtain the double generating
function~\eqref{eq:Stirling_DoubleGF} immediately from the
exponential formula (a set partition is a ``set of nonempty sets'' in
the language of species), and then by specializations, coefficient
extractions,~etc.\ recover all the previous results and even much
more.

\subsubsection{Coefficient extraction} Given a Z-statistic $\stat$, the first step we have
to perform in order to be able to apply Theorem~\ref{thm:GfunMean} is to try to find a
convenient expression for the power series $V_2^{\stat}(x)$.
  It turns out that for all Z-statistics which the author 
was able to locate in the literature,
the formal power series $V_2^{\stat}(x)$ is given by
\begin{align}
V_2^{\stat}(x)=P_0(x)+P_1(x) e^{x}+P_2(x)
e^{2x},\label{eq:canform_Vstat}
\end{align}
where the $P_i$  are polynomials. Therefore, we need to develop some tools
to cope with this situation.
 Elementary computations based on the formulas \eqref{eq:Bell_GF} and \eqref{eq:Stirling_VerticalGF} show that we have
\begin{align}
&e^{x}B(x)=B^{\prime}(x),\quad e^{2x}B(x)=B^{\prime\prime}(x)-B^{\prime}(x)\label{eq:expgf_der12_Bell}\\
&e^{x}S_{k-2}(x)=S_{k-2}(x)+(k-1)S_{k-1}(x)\label{eq:expgf_der1_Stirling}\\
&e^{2x}S_{k-2}(x)=S_{k-2}(x)+2(k-1)S_{k-1}(x)+k(k-1)S_{k}(x)\label{eq:expgf_der2_Stirling}.
\end{align}
This implies that, for $V_2^{\stat}(x)$ satisfying
\eqref{eq:canform_Vstat}, there exist  polynomials $Q_{i}$
and $R_{i}$, ${0\leq i\leq 2}$, such that
\begin{align}
V_2^{\stat}(x)B(x)&=Q_0(x) B(x)+Q_1(x) B^{\prime}(x)+Q_2(x) B^{\prime\prime}(x)\label{eq:canform_VstatB}\\
V_2^{\stat}(x)S_{k-2}(x)&=R_0(x) S_{k-2}(x)+R_1(x) S_{k-1}(x)+R_2(x)
S_{k}(x).\label{eq:canform_VstatS}
\end{align}
Suppose we are given a formal power series
$A(x)=\sum_{n\geq0}a_n\frac{x^{n}}{n!}$ and let $A^{(h)}(x)$ be the
 $h$-th formal derivative $A^{(h)}(x)$ of $A(x)$. Using the basic facts
\begin{align}
\left[\frac{x^{n}}{n!}\right]A^{(h)}(x)&=a_{n+h} \qquad
\left[\frac{x^{n}}{n!}\right]x^{i}A(x)=(n)_i a_{n-i},
\label{eq:extractionCoeff}
\end{align}
coefficient extraction from power series of the
form~\eqref{eq:canform_VstatB} or~\eqref{eq:canform_VstatS}
becomes a routine process. In particular, if $V_2^{\stat}(x)$
satisfies~\eqref{eq:canform_Vstat}, it follows
by Theorem~\ref{thm:GfunMean}
and~\eqref{eq:canform_VstatB}--\eqref{eq:extractionCoeff}
that $\mu_n$ and $\mu_{n,k}$ are given by
\begin{align}
\mu_n=\sum_{i\in
\mathbb{Z}}p_i(n)\frac{B_{n+i}}{B_{n}}\quad\text{and}\quad
\mu_{n,k}=\sum_{i,j\geq 0} q_{i,j}
(n)\frac{S_{n-i,k-j}}{S_{n,k}},\label{eq:canform_mean}
\end{align}
where the $p_i(n)$ and $q_{i,j} (n)$ are polynomials in $n$ which
are all zero except for a finite number of them.

\subsubsection{Asymptotics}
In order to obtain asymptotic approximations of the average values
$\mu_n$ and $\mu_{n,k}$, in view of~\eqref{eq:canform_mean}, 
it is important to consider quotients of Bell
numbers and Stirling numbers. We begin with Stirling numbers.

  It follows immediately from~\eqref{eq:Stirling_SummationFormula} that,
 as $n\rightarrow \infty$, we have
\begin{align}
S_{n,k}&=\frac{k^n}{k!}+O\left(\left(k-1\right)^n\right)
=\frac{k^n}{k!}\left(1+O\left(\left(1-\frac{1}{k}\right)^n\right)\right),\label{eq:Asymptotic_Stirling}
\end{align}
from which we deduce that, for $i,j\geq0$, we have
\begin{align*}
\frac{S_{n-i,k-j}}{S_{n,k}}&=\frac{(k-j)^{n-i}}{(k-j)!}\left(1+O\left(\left(1-\frac{1}{k-j}\right)^{n-i}\right)\right)
\frac{k!}{k^n}\left(1+O\left(\left(1-\frac{1}{k}\right)^n\right)\right)^{-1}\\
&=\frac{k!}{(k-j)!\left(k-j\right)^{i}}{\left(1-\frac{j}{k}\right)^{n}}
\left(1+O\left(\left(1-\frac{1}{k}\right)^n\right)\right).
\end{align*}
Distinct the cases $j=0$ and $j\geq 1$ leads to the following
result.
\begin{lem}\label{lem:Asymptotic_QuotientStirling}
For $i\geq0$ and $j\geq1$, we have, as $n\rightarrow\infty$,
\begin{align*}
\frac{S_{n-i,k}}{S_{n,k}}&=\frac{1}{k^i}+O\left(\left(1-\frac{1}{k}\right)^n\right)\\
\frac{S_{n-i,k-j}}{S_{n,k}}&=O\left(\left(1-\frac{j}{k}\right)^n\right).
\end{align*}
\end{lem}

 The asymptotics of Bell numbers is more delicate and is best performed by a saddle-point method.
It is probably the most famous application
of saddle-point techniques to combinatorial enumeration. The Bell
number $B_n$  satisfies (see \cite[Prop.~VIII.3]{FlSeg})
\begin{align*}
 B_n = n! \frac{e^{e^r-1}}{r^n\sqrt{2\pi r (r+1)e^r}}\left(1+O\left(e^{-r/5}\right)\right),
\end{align*}
where $r$ is defined implicitly by $re^r = n + 1$, so that $r =( \log n )\cdot (\log \log n) + o(1)$.
Instead of the latter formula, we  will use the following 
corollary given by Salvy and Shackell~\cite{Sal}:
\begin{align}
\frac{B_{n+1}}{B_{n}} &=\frac{n}{\log n} \left(1+\frac{\log\log
n}{\log n} (1+o(1))\right).\label{eq:Asymptotic_Bell_lem1}
\end{align}
Note that the latter identity was established because it gives an
asymptotic approximation of the average  number of blocks in a
random set partition of~$[n]$. It is also worth noting that there
exists a more precise asymptotic approximation (see~\cite{Sal})
which the reader can use to get more precise asymptotic
approximations for the average values computed in this paper.

\begin{lem}\label{lem:Asymptotic_QuotientBell}
For any integer $r\in \mathbb{Z}$, we have, as $n\rightarrow\infty$,
\begin{align}
\frac{B_{n+r}}{B_{n}} &=\left(\frac{n}{\log n}\right)^{r}
 \left (1+r\frac{\log\log n}{\log n} (1+o(1)) \right).\label{eq:Asymptotic_QuotientBell1}
\end{align}
\end{lem}

\pf
 Suppose we are given an integer $j$. It follows from
 \eqref{eq:Asymptotic_Bell_lem1} that
\begin{align}
  \frac{B_{n+j+1}}{B_{n+j}}&=\frac{n+j}{\log (n+j)}\left(1+\frac{\log\log (n+j)}{\log (n+j)} (1+o(1))
  \right).\label{eq:Formula_QuotientBell_shift}
\end{align}

By elementary calculus and the well-known formula $\log
(1+x)=x(1+o(1))$ as $x\rightarrow0$, it is easy to establish that
\begin{align}
\log (n+j)=\log n\left(1+O\left(\frac{1}{ n\log n}\right)\right)
\;\text{and}\;\; \frac{\log\log (n+j)}{\log (n+j)}=\frac{\log\log
n}{\log n}\,\left(1+o(1)\right),\label{eq:elementary_Asymptotic}
\end{align}
and thus
\begin{align}
\frac{n+j}{\log (n+j)}=\frac{n\left(1+\frac{j}{n}\right)}{\log
n}\left(1+O\left(\frac{1}{ n\log n}\right)\right)^{-1}=\frac{n}{\log
n}\left(1+O\left(\frac{1}{n}\right)\right).\label{eq:elementary_Asymptotic1}
\end{align}
Combining \eqref{eq:Formula_QuotientBell_shift},
\eqref{eq:elementary_Asymptotic}, and
\eqref{eq:elementary_Asymptotic1}, we obtain immediately
\begin{align*}
  \frac{B_{n+j+1}}{B_{n+j}}&=\frac{n}{\log n} \left(1+\frac{\log\log n}{\log n} (1+o(1)) \right),
\end{align*}
from which we deduce that, for any integer $r$, we have
\begin{align*}
  \frac{B_{n+r}}{B_{n}} &=\prod_{j=0}^{r-1}
  \frac{B_{n+j+1}}{B_{n+j}}
                         =\left(\left(\frac{n}{\log n}\right) \left(1+\frac{\log\log n}{\log n} (1+o(1)) \right)\right)^{r}\\
                        &=\left(\frac{n}{\log n}\right)^{r} \left (1+r\frac{\log\log n}{\log n} (1+o(1)) \right)
\end{align*}
since $(1+x)^{r}=1+r x + o\left(x\right)$ as $x\rightarrow0$. \qed

\subsection{A process}\label{sect:process}

 By Theorem~\ref{thm:GfunMean}, in order to get the average values $\mu_n$ and $\mu_{n,k}$ of a Z-statistic $\stat$,
we can proceed essentially as follows:
\begin{itemize}
\item[(i)] Compute $V_2^{\stat}(x)$.
\item[(ii)] Extract coefficients from the power series  $V_2^{\stat}(x)\,B(x)$ and $V_2^{\stat}(x)\,S_{k-2}(x)$.
\end{itemize}

 Actually, the above process seems to be nondeterministic, since it is a priori
 not clear whether $V_2^{\stat}(x)$ admits a simple expression (a finite expression involving elementary functions).
However, as we mentioned previously, for all Z-statistics on set partitions 
which the author was able to locate in the literature,
the formal power series $V_2^{\stat}(x)$ admits the simple form \eqref{eq:canform_Vstat}. 
For such statistics, one can
design a deterministic algorithm.\\

\textit{Input}: a Z-statistic $\stat$ such that  $V_2^{\stat}(x)=P_0(x)+P_1(x) e^{x}+P_2(x) e^{2x}$,
 where the $P_{i}$  are polynomials.

\textit{Output}: (a) the exact average values $\mu_n$ and $\mu_{n,k}$, written as finite expressions which involve
elementary functions, Bell and Stirling numbers; (b) asymptotic approximations of $\mu_n$ and $\mu_{n,k}$.\\

\textit{Step 1}. Compute explicitly  $V_2^{\stat}(x)$ (i.e., determine the polynomials $P_0$, $P_1$ and $P_2$).

\textit{Step 2}. Using~\eqref{eq:expgf_der12_Bell}--\eqref{eq:expgf_der2_Stirling},
write   $V_2^{\stat}(x)B(x)$ and $V_2^{\stat}(x)S_{k-2}(x)$ as
\begin{align*}
V_2^{\stat}(x)B(x)&=Q_0(x) B(x)+Q_1(x) B^{\prime}(x)+Q_2(x) B^{\prime\prime}(x)\\
V_2^{\stat}(x)S_{k-2}(x)&=R_0(x) S_{k-2}(x)+R_1(x) S_{k-1}(x)+R_2(x)
S_{k}(x),
\end{align*}
where the $Q_{i}$ and $R_{i}$ are polynomials.

\textit{Step 3}.   Using~\eqref{eq:extractionCoeff}, write the coefficients  of $V_2^{\stat}(x)B(x)$ and $V_2^{\stat}(x)S_{k-2}(x)$
as
\begin{align*}
 \left[\frac{x^{n}}{n!}\right]V_2^{\stat}(x)B(x)=\sum_{i\in \mathbb{Z}}p_i(n)B_{n+i}\;\;\text{and}\;\;
 \left[\frac{x^{n}}{n!}\right]V_2^{\stat}(x)S_{k-2}=\sum_{i,j\geq 0} q_{i,j} (n,k) S_{n-i,k-j},
\end{align*}
where the $p_i(n)$ and $q_{i,j} (n,k)$ are polynomials in $n$ and $k$ which
are all zero except for a finite number of them.

\textit{Step 4}. Divide the expressions obtained in the previous step
by $B_n$ and $S_{n,k}$ in order to get the average values $\mu_n$ and $\mu_{n,k}$,
which can be written in the form
\begin{align*}
 \left[\frac{x^{n}}{n!}\right]V_2^{\stat}(x)B(x)=\sum_{i\in \mathbb{Z}}p_i(n)\frac{B_{n+i}}{B_n}\;\;\text{and}\;\;
 \left[\frac{x^{n}}{n!}\right]V_2^{\stat}(x)S_{k-2}=\sum_{i,j\geq 0} q_{i,j} (n,k) \frac{S_{n-i,k-j}}{S_{n,k}}.
\end{align*}

\textit{Step 5}. Use Lemmas~\ref{lem:Asymptotic_QuotientStirling} and
\ref{lem:Asymptotic_QuotientBell} in the expressions obtained in the previous step
to get asymptotic approximations of  $\mu_n$ and $\mu_{n,k}$.\\

A proof that the above process is deterministic can be easily
extracted from the previous subsection (details are left to the
reader). To make things more concrete, we give an example in the next
subsection.

\subsection{Example: an inversion statistic on set partitions}\label{section:inversion}

 In his study of some $q$-Stirling numbers which arose from a $q$-exponential formula, Johnson~\cite{Johnson}
(see also~\cite{DeoSri}) introduced an inversion statistic on set partitions, denoted here
by $\inv$, which can be defined as follows. For $\pi\in\Pi_n$ with
standard form $\pi=B_1/B_2/\cdots\/B_k$, $\inv(\pi)$ is
the number of pairs $(i,j)$, $1\leq i<j\leq n$, such that $i$
belongs to a block to the right of the block containing $j$. Equivalently,
$\inv(\pi)$ is the number of inversions in the RGF $w(\pi)$ of $\pi$. Recall that an inversion in a word $w_1\ldots w_n$
is a pair $(i,j)$, $1\leq i<j\leq n$, such that $w_i>w_j$.  For instance, there are exactly 3
inversions in $\pi=14/25/3\equiv 1\,2\,3\,1\,2=w(\pi)$: the pairs $(2,4)$,
$(3,4)$ and $(3,5)$; hence $\inv(\pi)=3$. Note that 
$\inv=\occ_{2\,1}$ (where $\occ_{2\,1}$ is defined in Section~\ref{section:inversion}).\\

We propose to determine  exact and asymptotic average numbers of
inversions in a random set partition. Since  $\inv$ is  a
Z-statistic, we must first ``compute''
$V_2^{\inv}(x)$.

 By definition,  $V_2^{\inv}(x)$ is the exponential generating function
 of the sequence $\left(v_n\right)_{n\ge0}$
 with $v_n:=v^{\inv}_{n,2}=\sum_{\pi\in\Pi_n^2}\inv(\pi)$. For a combinatorialist,
it is usual to interpret $v_n$ as the number of pairs
$(\pi,a)$, where~$\pi$ is a 2-partition of $[n]$
and~$a=(i,j)$ is an inversion of~$\pi$. Such a pair will be
called an \emph{underlined} 2-partition.
 An underlined partition $(\pi,a)$
with $a=(i,j)$ can be identified with the restricted growth function
$w(\pi)$ in which the letters $w_i$, $w_j$ are colored. For
instance, there are 3 inversions in $\pi=14/25/3$: $a_1=(2,4)$,
$a_2=(3,4)$ and $a_3=(3,5)$. Then we have
$(\pi,a_1)\equiv 1\,\textbf{2}\,3\,\textbf{1}\,2$,
 $(\pi,a_2)\equiv 1\,2\,\textbf{3}\,\textbf{1}\,2$
and $(\pi,a_3)\equiv 1\,2\,\textbf{3}\,1\,\textbf{2}$.

 Then, using this correspondence, it is not hard
to see that an underlined 2-partition ${(\pi,a)}$ can be uniquely written as
\begin{align*}
(\pi,a)\equiv 1\,u_1 \, \2 \, u_2 \,\1\,u_3
\end{align*}
 with $u_1,u_2,u_3\in \{1,2\}^*$.  
Here and in the rest of the paper,  we write  $X^{*}$ for the free monoid generated by an alphabet $X$, which is
the set of words whose letters are in $X$. The empty word will be denoted $\epsilon$. 

Elementary counting (there are ${n-1 \choose 2}$ choices for the
positions of the underlined elements and $2^{n-3}$ choices for the
word $u_1 \,  u_2 \,u_3$) leads to
\begin{equation*}
                                v_n=\begin{cases} 
                                         {n-1 \choose 2}\,2^{n-3}, & \hbox{if $n\geq2$;} \\
                                         0, & \hbox{if $n\leq1$.}
                                       \end{cases}
\end{equation*}
 A straightforward computation, the details of which are left to the
 reader, then yields 
\begin{align}\label{eq:gfV_InvJohnson}
 V_2^{\inv}(x)&=\sum_{n\geq0}v_n\frac{x^n}{n!}=\left(\frac{x^2}{4}-\frac{x}{4}+\frac{1}{8}\right)e^{2x}-\frac{1}{8}.
\end{align}

 Since $V_2^{\inv}(x)$ can be written in the
 form~\eqref{eq:canform_Vstat}, 
we can continue the process
described in the previous subsection.\\

\textit{Step 2.} Using~\eqref{eq:expgf_der12_Bell} and~\eqref{eq:expgf_der2_Stirling}, it is straightforward to obtain
\begin{align*}
V_2^{\inv}(x)B(x)&=-\frac{1}{8} B(x)+\frac{1}{8} \left(-1+2 x-2 x^2\right) B^{\prime}(x)+\frac{1}{8}\left(1-2 x+2 x^2\right) B^{\prime\prime}(x)\\
V_2^{\inv}(x)S_{k-2}(x)&= \frac{1}{4} (-x+x^2)S_{k-2}(x)+\frac{1}{4} (k-1) \left(1-2 x+2 x^2\right)S_{k-1}(x)\\
&\qquad\qquad\hspace{0cm}+ \frac{1}{8}  k(k-1)\left(1-2 x+2
x^2\right) S_{k}(x).
\end{align*}

\textit{Step 3.} After routine coefficient extraction based on
\eqref{eq:extractionCoeff}, we obtain
\begin{align}
 \left[\frac{x^{n}}{n!}\right]V_2^{\inv}(x)B(x) &=\frac{1}{8}
B_{n+2}-\left(\frac{n}{4}+\frac{1}{8}\right)B_{n+1}
+\left(\frac{n^2}{4}-\frac{1}{8}\right)
B_{n}-\frac{(n)_2}{4}B_{n-1} \label{eq:V*B_inv} \\[0.3cm]
 \left[\frac{x^{n}}{n!}\right]V_2^{\inv}(x)S_{k-2}(x)
&= \frac{1}{4} \left(-n S_{n-1,k-2}+(n)_2 S_{n-2,k-2}  \right) \nonumber\\
&\hspace{1cm} +\frac{1}{4} (k-1) \left(S_{n,k-1}-2 n S_{n-1,k-1}+ 2 (n)_2 S_{n-2,k-1} \right)\nonumber\\
&\hspace{1cm}+ \frac{1}{8}  k(k-1)\left(S_{n,k}-2 n S_{n-1,k}+ 2 (n)_2
S_{n-2,k} \right).\nonumber
\end{align}

It is possible to simplify the expression obtained for $\left[\frac{x^{n}}{n!}\right]V_2^{\inv}(x)S_{k-2}(x)$. First, replace each occurrence
of the left hand sides of the three identities
\begin{align}
S_{n-2,k-2}&=S_{n-1,k-1}-(k-1)S_{n-2,k-1}\label{eq:StirlingSimplification1}\\
S_{n-1,k-2}&= S_{n,k-1}-(k-1)S_{n-1,k-1}\label{eq:StirlingSimplification2}\\
S_{n-2,k}&=\frac{1}{k}\left(S_{n-1,k}-S_{n-2,k-1}\right)\label{eq:StirlingSimplification3}
\end{align}
by the  corresponding right hand sides. 
Then, in the identity obtained from the previous manipulation, replace each occurrence
of the left hand side of the identity
\begin{align}
S_{n-1,k}&=\frac{1}{k}\left(S_{n,k}-S_{n-1,k-1}\right)\label{eq:StirlingSimplification4}
\end{align}
by its right hand side. 
This gives
\begin{align}
 \left[\frac{x^{n}}{n!}\right]V_2^{\inv}(x)S_{k-2}(x)
 &=\left(\frac{1}{4} n \left(n-k-1\right) \left(1-\frac{1}{k}\right)+\frac{1}{8}  k(k-1)\right)S_{n,k} \label{eq:V*S_inv}\\
&\qquad -\frac{1}{4}\left(n+1-k
\right)S_{n,k-1}+\frac{n(n-1)}{4k}S_{n-1,k-1}.\nonumber
\end{align}
 Note that identities \eqref{eq:StirlingSimplification1}--\eqref{eq:StirlingSimplification4} could be obtained in an elementary way
from~\eqref{eq:Stirling_Recurrence}.\\

\textit{Steps 4 and 5.}  Divide the expression \eqref{eq:V*B_inv} by
$B_n$ to obtain that \emph{the average number of inversions $\mu_n$ in
a random set partition of~$[n]$ is given by}
\begin{align*}
\mu_n&=\frac{1}{8}\frac{B_{n+2}}{B_{n}}
-\left(\frac{n}{4}+\frac{1}{8}\right)\frac{B_{n+1}}{B_{n}}
+\left(\frac{n^2}{4}-\frac{1}{8}\right)-\frac{(n)_2}{4}\frac{B_{n-1}}{B_{n}}.
\end{align*}
 Using~Lemma~\ref{lem:Asymptotic_QuotientBell}, it is straightforward
 to arrive at the asymptotic approximation
\begin{align*}
\mu_{n}&=\frac{n^2}{4}\left( 1-\frac{1}{\log n}-\frac{\log\log
n}{\left(\log n\right)^2}\left(1+o\left(1\right)\right) \right),\quad   n\rightarrow\infty.
\end{align*}

 Similarly, divide expression \eqref{eq:V*S_inv} by  $S_{n,k}$ to obtain that
\emph{the average number of inversions $\mu_{n,k}$  in a
random set partition of~$[n]$ into~$k$ blocks is given by}
\begin{align*}
\mu_{n,k}&=\frac{1}{4} n \left(n-k-1\right)
\left(1-\frac{1}{k}\right)+\frac{1}{8}  k(k-1)
-\frac{1}{4}\left(n+1-k
\right)\frac{S_{n,k-1}}{S_{n,k}}\\
&\qquad +\frac{n(n-1)}{4k}\frac{S_{n-1,k-1}}{S_{n,k}}.
\end{align*}
  Use of~Lemma~\ref{lem:Asymptotic_QuotientStirling} then produces the
  asymptotic approximation
\begin{align*}
\mu_{n,k} &=\frac{1}{4} n \left(n-k-1\right)
\left(1-\frac{1}{k}\right)+\frac{1}{8}  k(k-1)+ O \left(
n\left(1-\frac{1}{k}\right)^n \right)\\
&=\frac{1}{4} n \left(n-k-1\right)
\left(1-\frac{1}{k}\right)+\frac{1}{8}  k(k-1)+ o(1),
\quad \text{as $n\rightarrow\infty$}.
\end{align*}

\subsection{Proof of Theorem~\ref{thm:GfunMean}}

 Let $\stat$ be a Z-statistic defined on set partitions.
It is easy to see that~\eqref{eq:ThmGfunMean} is a corollary
of~\eqref{eq:ThmGfunMeanBlock}. Indeed, by definition of the average
value, we have
\begin{align}\label{eq:Meanblock_Mean}
 \mu_n=\frac{ \sum_{\pi\in \Pi_n} \stat(\pi)}{B_n}=\frac{
\sum_{k\geq 1}\sum_{\pi\in \Pi_n^k}
\stat(\pi)}{B_n}=\frac{1}{B_n}\sum_{k\geq 1} S_{n,k}\, \mu_{n,k}.
\end{align}
Suppose that~\eqref{eq:ThmGfunMeanBlock} is true. Then
we have
\begin{align*}
 \mu_n&=\frac{1}{B_n}\sum_{k\geq 1} S_{n,k}\, \mu_{n,k}
      =\frac{1}{B_n}\sum_{k\geq 1}  \left[\frac{x^{n}}{n!}\right]V_2^{\stat}(x)S_{k-2}(x)\\
       &=\frac{1}{B_n}  \left[\frac{x^{n}}{n!}\right]V_2^{\stat}(x)\left(\sum_{k\geq 1} S_{k-2}(x)\right)
     =\frac{1}{B_n}  \left[\frac{x^{n}}{n!}\right] V_2^{\stat}(x) B(x).
\end{align*}
Therefore, in order to prove Theorem~\ref{thm:GfunMean}, we just
have to prove~\eqref{eq:ThmGfunMeanBlock}. Set
$v_m:=v_{m,2}^{\stat}$. By the definitions of average value and
product of two formal power series,
Equation~\eqref{eq:ThmGfunMeanBlock} is equivalent to 
\begin{align}
\sum_{\pi\in \Pi_n^k} \stat(\pi)&=\sum_{m=0}^n {n\choose
m}S_{n-m,k-2}\,v_m.\label{eq:Ank}
\end{align}
 We will prove the latter identity.
  Let $DISJ$ be the set of pairs of nonempty and
disjoints subsets of $\mathbb{P}$, i.e.,
$$
DISJ=\{\{A,B\}\,/\,A,B\subseteq \mathbb{P},\;A,B\neq\emptyset,\;
A\cap B=\emptyset\}.
$$
 Then, by definition of a Z-statistic, we have
\begin{align}
\sum_{\pi\in \Pi_n^k} \stat(\pi)&=\sum_{\pi\in \Pi_n^k} \sum_{A,B\in
\pi}\stat(\st(\{A,B\})) =\sum_{A,B\subseteq [n]} \sum\limits_{\substack{\pi\in \Pi_n^k \\
A,B\in \pi}} \stat(\st(\{A,B\}))\nonumber\\
&=\sum\limits_{\substack{A,B\subseteq [n] \\
\{A,B\}\in DISJ}} \stat(\st(\{A,B\}))\,S_{n-|A\cup B|,k-2},\label{eq:MainThm1}
\end{align}
where the last equality follows from the fact that, for $A,B\subseteq
[n]$, the number~$p_{n,k}(A,B)$ of partitions~$\pi$ in~$\Pi_n^k$ such
that $A,B\in\pi$ is given by $$p_{n,k}(A,B)=\begin{cases}
                               S_{n-|A\cup B|,k-2}, & \hbox{if $\{A,B\}\in DISJ$;} \\
                               0, & \hbox{otherwise.}
\end{cases}
$$
Now, suppose we are given a set $C\subseteq \mathbb{P}$ with $|C|=m$. We have
 \begin{align}
\sum\limits_{\substack{\{A,B\}\in DISJ\\
A\cup
B=C}}\stat(\st(\{A,B\}))=\sum_{\tau\in\Pi^2(C)}\stat(\st_C(\tau))=\sum_{\pi\in\Pi_m^2}\stat(\pi)=v_m,\label{eq:MainThm2}
\end{align}
where the first equality follows from the definition of $DISJ$ and the second
is a consequence of the fact that $\st_C$ send $\Pi^2(C)$ bijectively onto $\Pi_m^2$.

Combining \eqref{eq:MainThm1} and \eqref{eq:MainThm2}, we obtain
\begin{align*}
\sum_{\pi\in \Pi_n^k} \stat(\pi)
&=\sum\limits_{\substack{A,B\subseteq [n] \\
\{A,B\}\in DISJ}} \stat(\st(\{A,B\}))\,S_{n-|A\cup B|,k-2}\\
&=\sum_{C\subseteq[n]} \sum\limits_{\substack{\{A,B\}\in DISJ\\
A\cup B=C}}\stat(\st(\{A,B\}))\,S_{n-|A\cup B|,k-2}\\
&=\sum_{C\subseteq[n]} v_{|C|}\,S_{n-|C|,k-2},
\end{align*}
which is obviously equivalent to~\eqref{eq:Ank}. This concludes the proof.

\begin{rmk}
Using the language of what is called the theory of species
one can give a quick proof
of~Theorem~\ref{thm:GfunMean}. We have preferred to avoid this
terminology in order to keep the paper self-contained.
\end{rmk}


\section{Carlitz's $q$-Stirling distributions}

\subsection{Proof of Theorem~\ref{thm:Mean_Qstirling}}

  We will work here with the Z-statistic $\los$,
which has the $q$-Stirling distribution $S(q)$ (see~\cite{Mi,WW}). 
Recall that $\los$ can be defined for a set partition
$\pi=B_1/B_2/\ldots/B_k$ in standard form by
$$\los(\pi)=|B_2|+2 |B_3|+\cdots+(k-1)|B_k|.$$

We have to start by ``computing''
$V_2^{\los}(x)$.\\

\emph{Computation of $V_2^{\los}(x)$.} By definition,
$V_2^{\los}(x)$ is the exponential generating function of the
sequence $\left(v_n\right)_{n\ge0}$ with $v_n:=v^{\los}_{n,2}=\sum_{\pi\in\Pi_n^2}\los(\pi)$.

 For a 2-partition $\pi=B_1/B_2$ written in standard form,
we have $\los(\pi)=|B_2|$.  Thus, we can see  $v_n$ as the number of
pairs $(\pi,i)$, where~$\pi=B_1/B_2$  is a 2-partition of $[n]$
(written in standard form) and~$i\in B_2$. The number of such pairs
is easily shown to be $(n-1)2^{n-2}$. Indeed, we have $n-1$ choices
 for the element $i$ (any of the elements of $\{2,3,\ldots, n\}$)
 and then $2^{n-2}$ choices for the set $B_2\setminus\{i\}$
 (any of the subsets of $\{2,3,\ldots, n\}\setminus\{i\}$).
Therefore, we have
$$
v_n=\sum_{\pi\in\Pi_n^2}\los(\pi)=\begin{cases} 
                                         0, & \hbox{if $n\leq1$;} \\
                                         (n-1)\,2^{n-2}, & \hbox{if $n\geq2$.}
                                     \end{cases}
$$
 A straightforward computation leads to
\begin{align}
 V_2^{\los}(x)&=\sum_{n\geq0}v_n\frac{x^n}{n!}=\sum_{n\geq2}(n-1)2^{n-2}\frac{x^n}{n!}
=\frac{1}{4}+\left(-\frac{1}{4}+\frac{x}{2}\right)e^{2x}.\label{eq:gfV_LOS}
\end{align}

 Since $V_2^{\los}(x)$ can be written in the
 form~\eqref{eq:canform_Vstat}, 
we can continue the process
described in Section~\ref{sect:process}.\\

\textit{Step 2.} Using~\eqref{eq:expgf_der12_Bell}  and
\eqref{eq:expgf_der2_Stirling}, we can write
\begin{align*}
V_2^{\los}(x)B(x)
&=\frac{1}{4} B(x) +\left(-\frac{x}{2}+\frac{1}{4}\right) B^{\prime}(x) +\left(\frac{x}{2}-\frac{1}{4}\right)B^{\prime\prime}(x)\\
V_2^{\los}(x)S_{k-2}(x) &=\frac{x}{2}S_{k-2}(x)+\frac{1}{2}(k-1)
\left(-1+2 x\right) S_{k-1}(x)+\frac{1}{4} k(k-1) \left(-1+2
x\right) S_{k}(x).
\end{align*}

\textit{Step 3.} After routine coefficient extraction based on
\eqref{eq:extractionCoeff}, we obtain
\begin{align}
 \left[\frac{x^{n}}{n!}\right]V_2^{\los}(x)B(x) &=-\frac{1}{4}
B_{n+2}+\left(\frac{n}{2}+\frac{1}{4}\right)B_{n+1}
+\left(-\frac{n}{2}+\frac{1}{4}\right) B_{n}\label{eq:V*B_los}\\
 \left[\frac{x^{n}}{n!}\right]V_2^{\los}(x)S_{k-2}(x)
&=\frac{n}{2}S_{n-1,k-2} +\frac{1}{2}(k-1) \left(-S_{n,k-1}+2 n
S_{n-1,k-1}\right)\nonumber\\
&\hspace{1cm}+\frac{1}{4} k(k-1) \left(-S_{n,k}+2 n S_{n-1,k}\right).\nonumber
\end{align}

It is possible to simplify the expression obtained
 for $\left[\frac{x^{n}}{n!}\right]V_2^{\los}(x)S_{k-2}(x)$ by replacing  each occurrence
of the left hand sides of identities~\eqref{eq:StirlingSimplification2} and~\eqref{eq:StirlingSimplification4}
by the  corresponding right hand sides.
This gives
\begin{align}
 \left[\frac{x^{n}}{n!}\right]V_2^{\los}(x)S_{k-2}(x)
&=\left(\frac{1}{2} n(k - 1)- \frac{1}{4}k(k - 1)\right)S_{n,k}
  + \frac{1}{2}(n+1-k)S_{n,k-1}. \label{eq:V*S_los}
\end{align}

\textit{Steps 4 and 5.} Division of expression~\eqref{eq:V*B_los}
by $B_n$ gives the exact
value~\eqref{eq:ExactMean_qstirling}  of $\mu_n$, while its
asymptotic approximation~\eqref{eq:AsymptoticMean_qstirling} is
easily obtained from Lemma~\ref{lem:Asymptotic_QuotientBell}.

 Similarly, division of expression~\eqref{eq:V*S_los}
by $S_{n,k}$ gives the exact
value~\eqref{eq:ExactMeanBlock_qstirling} of $\mu_{n,k}$, while,
after a routine computation based on
Lemma~\ref{lem:Asymptotic_QuotientStirling}, we obtain the asymptotic
approximation
\begin{align*}
\mu_{n,k}&=\frac{1}{2} n(k - 1)- \frac{1}{4}k(k - 1)+ O \left(
n\left(1-\frac{1}{k}\right)^n \right),\quad n\rightarrow\infty,
\end{align*}
which is a refinement of~\eqref{eq:AsymptoticMeanBlock_qstirling}.
This concludes the proof of Theorem~\ref{thm:Mean_Qstirling}.

\subsection{Some remarks}
There is an ``explicit'' expression for the generating function of the $q$-Stirling
distribution $S(q)$ (see e.g.~\cite{Gou}), namely,
\begin{align*}
 S_{n,k}(q)=\frac{1}{[k]_q!}\sum_{j=1}^{k} (-1)^{k-j} {k\brack j}_q  q^{\binom{k-j}{2}}[j]_q^n,
\end{align*}
where  ${n\brack k}_q$ is the usual $q$-binomial coefficient. 
It seems difficult to recover the exact average values
of the $q$-Stirling distributions from the above expression.
However, it is well-known that ''generating functions find averages`` (see e.g.~\cite{Wilf}). 
In the case of our $q$-Stirling numbers, it is easy to derive from the
recurrence~\eqref{eq:stirling1} that
\begin{align}\label{eq:GenFunc_qStirling}
F_k(x,q):=\sum_{n\geq 0}S_{n,k}(q)x^n=\frac{q^{k \choose 2}\,x^k}{(1-x)(1-[2]_q
x)\cdots(1-[k]_q x)}.
\end{align}
 It is possible to recover the results presented in this section from
 this generating function identity. We can even
 obtain many more results, such as values of variances and 
some limit laws (see~\cite{KaFutur}).


\section{Number of linear crossings}

\subsection{Proof of Theorem~\ref{thm:Mean_CrossingLinear}}\label{sect:defTk}

\emph{Computation of $V_2^{\crol}(x)$.} 
By definition, $V_2^{\crol}(x)$ is the exponential generating function of the
sequence $\left(v_n\right)_{n\ge0}$ with
$v_n:=v^{\crol}_{n,2}=\sum_{\pi\in\Pi_n^2}\crol(\pi)$.

  We can interpret $v_n$ as the number of pairs $(\pi,C)$,
where~$\pi$ is a 2-partition of $[n]$ and~$C$ is a linear crossing
of~$\pi$. Such a pair will be called an \emph{underlined}
2-partition.

 An underlined partition $(\pi,C)$ with $C=(i,j)(k,\ell)$ can be identified
with the restricted growth function $w(\pi)$ in which the letters $w_i$, $w_j$, $w_k$ and $w_{\ell}$
are colored.  For instance, there are 3 linear crossings in
$\pi=1\,4/2\,5/3\,6\,7$: ${C_1}=(1,4)(2,5)$,
 ${C_2}=(1,4)(3,6)$ and ${C_3}=(2,5)(3,6)$. Then we have
$(\pi,{C_1})\equiv\textbf{1}\,\textbf{2}\,3\,\textbf{1}\,\textbf{2}\,3\,3$,
 $(\pi,{C_2})\equiv\textbf{1}\,2\,\textbf{3}\,\textbf{1}\,2\,\textbf{3}\,3$
and
${(\pi,{C_3})\equiv1\,\textbf{2}\,\textbf{3}\,1\,\textbf{2}\,\textbf{3}\,3}$.

 Using this correspondence, it is not hard
to see that an underlined 2-partition~$(\pi,C)$ can be decomposed uniquely in
one of the following two ways:
\begin{itemize}
\item[(a)]$(\pi,C)\equiv u_1 \,\1 \,2^{k_1} \,\2\,\1\, 1^{k_2}\, \2\, u_2$,
\item[(b)]$(\pi,C)\equiv v_1\, \2 \,1^{\ell_1}\, \1\,\2 \,2^{\ell_2}\, \1\, v_2$,
\end{itemize}
 with $u_1=\epsilon$ or $u_1\in 1\{1,2\}^*$, $u_2\in \{1,2\}^*$ and $k_1,k_2\geq
 0$;
 $v_1\in 1\{1,2\}^*$, $v_2\in \{1,2\}^*$ and $\ell_1,\ell_2\geq
 0$.

 Simple rules of counting and manipulations of generating functions (see e.g.~\cite{FlSeg,Wilf}) then show
that the ordinary generating functions $F^{(a)}(x)$ and $F^{(b)}(x)$
for underlined\linebreak {2-partitions} of $[n]$, $n\in\mathbb{N}$ admitting a decomposition of
type~$(a)$ and $(b)$, respectively, are given by
\begin{align*}
F^{(a)}(x)&=\left(1+x\,\frac{1}{1-2x}\right)\,x\,\left(\frac{1}{1-x}\right)\,x^2\,\left(\frac{1}{1-x}\right)\,x\,\left(\frac{1}{1-2x}\right)
=\frac{x^4(1-x)}{(1-x)^2(1-2x)^2}\\
F^{(b)}(x)&=\left(x\,\frac{1}{1-2x}\right)\,x\,\left(\frac{1}{1-x}\right)\,x^2\,
\left(\frac{1}{1-x}\right)\,x\,\left(\frac{1}{1-2x}\right)
=\frac{x^5}{(1-x)^2(1-2x)^2}.
\end{align*}
Thus, the number of underlined 2-partitions of $[n]$, i.e., $v_n$, is
just the coefficient of $x^{n}$ in the power series
\begin{align*}
F(x)=F^{(a)}(x)+F^{(b)}(x)=\frac{x^4}{(1-x)^2(1-2x)^2}.
\end{align*}
A routine coefficient extraction based on the partial fraction
decomposition
\begin{align*}
F(x)=\frac{x^4}{(1-x)^2(1-2x)^2}=x^{4}\left(\frac{1}{(1 -
x)^2}+\frac{4}{1 - x}+\frac{4}{(1 -2 x)^2} -\frac{8}{1 -2 x}\right)
\end{align*}
 gives
$$
v_n=\sum_{\pi\in\Pi_n^2}\crol(\pi)=\begin{cases} 
                                         0, & \hbox{if $n\leq3$;} \\
                                         (n-5)\,2^{n-2}+ n + 1, & \hbox{if $n\geq4$.}
                                     \end{cases}
$$

 A straightforward computation then leads to
\begin{align}\label{eq:gfV_CrossingLinear}
 V_2^{\crol}(x)&=\sum_{n\geq0}v_n\frac{x^n}{n!}=\left(\frac{x}{2}-\frac{5}{4}\right)e^{2x}+(x+1)\,e^x+\frac{1}{4}.
\end{align}


 Of course, there are many other ways to compute the coefficients $v_n$.
We now present  another way based on a decomposition of
2-partitions which is particularly well-adapted to the computation
of the coefficients $v_n$ in the case of circular
crossings.

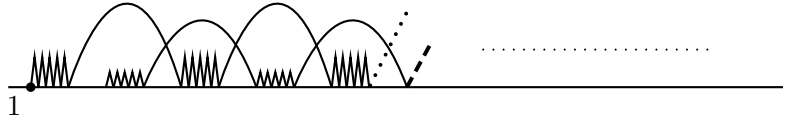
\begin{figure}
\begin{center}
\begin{pspicture}(0,0)(10,1.6) 
 \def\zigzagA{\multirput(0,0)(.1,0){5}{\psline(0,0)(0.05,.4)(.1,0)}}
 \def\zigzagB{\multirput(0,0)(.1,0){5}{\psline(0,0)(0.05,.2)(.1,0)}}
 \def\clocheA{\psbezier(0.5,0)(0.5,0)(1.25,2.5)(2,0)}
 \def\clocheB{\psbezier(0.5,0)(0.5,0)(1.25,2)(2,0)}
 \multirput(0,0)(2,0){2}{\clocheA}
 \multirput(0,0)(2,0){3}{\zigzagA}
 \multirput(1,0)(2,0){2}{\zigzagB}
 \multirput(1,0)(2,0){2}{\clocheB}
  \rput(4.5,0){\psline[linestyle=dotted,linewidth=1.5pt](0,0)(0.5,1)}
 \rput(5,0){\psline[linestyle=dashed,linewidth=1.5pt](0,0)(0.3,0.55)}
\psline(-0.3,0)(10,0) \psline[linestyle=dotted](6,0.5)(9,0.5)
 \uput[dl](0,0){\small $1$}
 \psdots(0,0)
\end{pspicture}
\caption{Linear representation of a 2-partition}
\label{fig:Linear2partition}
\end{center}
\end{figure}

 Via its restricted growth function, a 2-partition can be viewed
 as a word on the alphabet~$\{1,2\}$ whose first letter is 1. It is then immediate
to see that each 2-partition $\pi$ can be decomposed uniquely as
\begin{align*}
\pi&=1^{\ell_1} 2^{m_1} 1^{\ell_2} 2^{m_2}\cdots1^{\ell_k}
2^{m_k}1^{\ell_{k+1}}
\end{align*}
 for some $k\geq1$,  where $\ell_j,m_j\geq1$ for  $1\leq j\leq k$, and $\ell_{k+1}\geq0$.

For $k\geq1$, denote by $\T_k$ the set of 2-partitions with such a
decomposition. A moment's thought (see
Figure~\ref{fig:Linear2partition}) will convince the reader that
\begin{itemize}
\item if $\pi\in\T_1$, $\crol(\pi)=0$,
\item if $\pi\in\T_k$, $k\geq2$, $\crol(\pi)=\begin{cases} 
                                                     2 k-3, & \hbox{if $\ell_{k+1}=0$;} \\
                                                     2 k-2, & \hbox{if $\ell_{k+1}\geq1$.}
                                                \end{cases}
$
\end{itemize}

 It is easy to see from the definition of $\T_k$ that the ordinary generating functions
 $G_k^{(a)}(x)$ and $G_k^{(b)}(x)$ of 2-partitions of $[n]$, $n\in\mathbb{N}$,  in $\T_k$
 satisfying  $\ell_{k+1}=0$ and $\ell_{k+1}\geq1$, respectively, are given by
$G_k^{(a)}(x)=\left(\frac{x}{1-x}\right)^{2k}$ and
$G_k^{(b)}(x)=\left(\frac{x}{1-x}\right)^{2k+1}$.  It follows that for
$k\geq2$ we have  
\begin{align*}
\sum_{\pi \in\T_k} \cro(\pi)\,x^{\|\pi\|}
&=(2k-3)G_k^{(a)}(x)+(2k-2)G_k^{(b)}(x)\\
&=(2k-3)\,\left(\frac{x}{1-x}\right)^{2k}+(2k-2)\,\left(\frac{x}{1-x}\right)^{2k+1},
\end{align*}
where we have set $\|\pi\|=n$ for $\pi\in\Pi_n$.

Using the latter equality and the formal identity
$\sum_{k\geq1}k\,y^{k}=\frac{y}{(1-y)^2}$, after a
routine computation we obtain
\begin{align*}
F(x)=\sum_{n\geq 0}v_n\,x^n&=\sum_{k\geq2}\sum_{\pi \in\T_k}
\crol(\pi)\,x^{\|\pi\|}=\frac{x^4}{(1-x)^2(1-2x)^2},
\end{align*}
as it was previously derived.\\

Since $V_2^{\crol}(x)$ (see \eqref{eq:gfV_CrossingLinear}) can be
written in the form~\eqref{eq:canform_Vstat}, 
we can continue the process
described in Section~\ref{sect:process}.\\

\textit{Step 2.} Using~\eqref{eq:expgf_der12_Bell},
\eqref{eq:expgf_der1_Stirling} and \eqref{eq:expgf_der2_Stirling},
we can write
\begin{align*}
V_2^{\crol}(x)B(x)
&=\frac{1}{4}B(x)+\frac{1}{4} (9+2 x)B^{\prime}(x)+\frac{1}{4}(-5+2 x) B^{\prime\prime}(x)\\
V_2^{\crol}(x)S_{k-2}(x)&= \frac{3 x}{2}S_{k-2}(x)+\frac{1}{2} (k-1)
(-3+4 x)S_{k-1}(x)+\frac{1}{4} k(k-1) (-5+2 x)S_{k}(x).
\end{align*}

\textit{Step 3.} After routine coefficient extraction based on
\eqref{eq:extractionCoeff}, we obtain
\begin{align}
\left[\frac{x^{n}}{n!}\right]V_2^{\crol}(x)B(x) &=-\frac{5}{4}
B_{n+2}+\left(\frac{n}{2}+\frac{9}{4}\right)B_{n+1}
+\left(\frac{n}{2}+\frac{1}{4}\right) B_{n}\label{eq:V*B_crol}\\
\left[\frac{x^{n}}{n!}\right]V_2^{\crol}(x)S_{k-2}(x) &=\frac{3
n}{2}S_{n-1,k-2}+\frac{1}{2} (k-1) \left(-3 S_{n,k-1}+4 n
S_{n-1,k-1}\right)\nonumber\\
&\hspace{1cm}+\frac{1}{4} k(k-1) \left(-5 S_{n,k}+2 nS_{n-1,k}\right).\nonumber
\end{align}
It is possible to simplify the expression obtained
 for $\left[\frac{x^{n}}{n!}\right]V_2^{\crol}(x)S_{k-2}(x)$ by replacing each occurrence
of the left hand sides of identities~\eqref{eq:StirlingSimplification2} and~\eqref{eq:StirlingSimplification4}
by the corresponding right hand sides.
This gives
\begin{align}
 \left[\frac{x^{n}}{n!}\right]V_2^{\crol}(x)S_{k-2}(x)
&=\left(\frac{1}{2} n(k - 1)- \frac{5}{4}k(k - 1)\right)S_{n,k}
  + \frac{3}{2}(n+1-k)S_{n,k-1}.\label{eq:V*S_crol}
\end{align}

\textit{Steps 4 and 5.} Division of expression~\eqref{eq:V*B_crol}
by $B_n$ gives the exact value~\eqref{eq:ExactMean_crol}
of $\mu_n$, while its asymptotic
approximation~\eqref{eq:AsymptoticMean_crol} is easily obtained
from Lemma~\ref{lem:Asymptotic_QuotientBell}.

 Similarly, division of expression~\eqref{eq:V*S_crol}
by $S_{n,k}$ gives the exact
value~\eqref{eq:ExactMeanBlock_crol} of $\mu_{n,k}$, while, after a
routine computation based on
Lemma~\ref{lem:Asymptotic_QuotientStirling}, we get the asymptotic
approximation
\begin{align*}
\mu_{n,k}&=\frac{1}{2}n(k-1)-\frac{5}{2}{k\choose 2}+ O \left(
n\left(1-\frac{1}{k}\right)^n \right),\quad n\rightarrow\infty,
\end{align*}
which is a refinement of~\eqref{eq:AsymptoticMeanBlock_crol}. This
concludes the proof of Theorem~\ref{thm:Mean_CrossingLinear}.

\subsection{Some Remarks}

It is worth noting that there has been considerable interest in
studying crossings in matchings (set partitions each block of which has
exactly two elements) and set partitions.

In the case of the number of crossings in matchings, a remarkable
formula, often called the Touchard--Riordan formula, was made explicit by
Riordan~\cite{Ri}, who also mentioned the exact values of the average and the variance
of the number of crossings in matchings. Though they are ``contained''
in the Touchard--Riordan formula, it is not a priori clear whether
they admit a simple form. Later, Flajolet and Noy~\cite{FlN} used a
certain decomposition of the number of crossings to give a ``direct''
proof for the average value. This decomposition was later
generalized by the author~\cite{Ka} to compute the average number of linear 
crossings in a random set partition and is at the origin of this work. 

Let $T_{n,k}(q)$ be the generating function of set partitions
of~$[n]$ into~$k$~blocks with respect to the number of linear crossings,
i.e., $T_{n,k}(q)=\sum_{\pi\in\Pi_n^k} q^{\crol(\pi)}$. 
Biane~\cite{Bi} (see also~\cite{KaZe}) found a continued fraction expansion for the 
generating function $\sum_{n\geq k\geq 0} T_{n,k}(q)\, a^kt^n$.  
Recently, Stanton, Zeng and the author~(see~Equation~28 in~\cite{KaStZe}) proved that  
\begin{align}\label{eq:gf_crol}
\sum_{n\geq k\geq 0} T_{n,k}(q)\,a^k\,t^n&=\sum_{k=0}^\infty
\frac{(aqt)^k}{\prod_{i=1}^k(q^i-q^i[i]_qt+a(1-q)[i]_qt)},
\end{align}
from which they derive the following remarkable
formula for $T_{n,k}(q)$ (see~Equation~32 in~\cite{KaStZe}):
\begin{align*}
T_{n,k}(q)
&=\sum_{j=1}^{k} (-1)^{k-j} \frac{[j]_q^n}{[j]_q!} \sum_{i=0}^{k-j} \frac{(1-q)^{i}}{[k-j-i]_{q}!}
 q^{\binom{k-j-i+1}{2}-kj}\biggl( \binom{n}{i} q^{j}+\binom{n}{i-1}\biggr).
\end{align*}
Another remarkable formula for $T_{n,k}(q)$ was established 
recently by Josuat-Verg\`es and Rubey~\cite{JosRu}:
{\small
\begin{align*}
 T_{n,k}(q)= \frac{1}{(1-q)^{n-k}}\sum_{j=0}^{k}\sum_{i=j}^{n-k}(-1)^{i}
\left( {n\choose k+i}  {n\choose k-j} - {n\choose k+i+1}  {n\choose
k-j-1}\right){i\brack j}_q q^{{j+1 \choose 2}}.
\end{align*}}
Again, though the average number of linear crossings is ``contained''
in the above two formulas, it is not a priori clear whether
they admit a simple form. However, it is possible  to recover the 
exact average values presented in this section
and even exact values for the variance of the number of linear
crossings (see~\cite{KaFutur}) from the generating function~\eqref{eq:gf_crol}. 
Moreover, a limit law will be given in~\cite{KaFutur}.

Finally, we want to point out that there exist
combinatorial parameters which have the same distribution (hence,
the same average value) as the number of linear crossings on each~$\Pi_n^k$. 
This is the case  of the number of nestings of two
arcs~\cite{KaZe} and the major index for set partitions introduced
in~\cite{ChenGe}.


\section{Number of circular crossings}

\subsection{Proof of Theorem~\ref{thm:Mean_CrossingCircular}}

\emph{Computation of $V_2^{\croc}(x)$.} By definition,
$V_2^{\croc}(x)$ is the exponential generating function of the
sequence $\left(v_n\right)_{n\ge0}$ with
$v_n:=v^{\croc}_{n,2}=\sum_{\pi\in\Pi_n^2}\croc(\pi)$.

   Suppose  we are given a 2-partition $\pi$ in $\T_k$.
By definition of $\T_k$ (see Section~\ref{sect:defTk}), this means that there exist integers
$\ell_j,m_j\geq1$ for  $1\leq j\leq k$, and $\ell_{k+1}\geq0$ such
that
\begin{align*}
\pi&=1^{\ell_1} 2^{m_1} 1^{\ell_2} 2^{m_2}\cdots1^{\ell_k}
2^{m_k}1^{\ell_{k+1}}.
\end{align*}

\begin{figure}[h!]
\begin{center}
\begin{pspicture}(-1.7,-1.7)(1.7,1.7)
\scalebox{0.85}{ \psdots(0,2)
 \uput[ur](0,2){\scriptsize 1}
\SpecialCoor
 \pscircle(0,0){2}
\psarc[linestyle=dotted](0,0){1.5}{160}{-50}
\def\zigzagCa{\psline(2;0)(1.75;2)(2;4)(1.75;6)(2;8)(1.75;10)(2;12)(1.75;14)(2;16)(1.75;18)(2;20)}
\def\zigzagCb{\psline(2;0)(1.87;2)(2;4)(1.87;6)(2;8)(1.87;10)(2;12)(1.87;14)(2;16)(1.87;18)(2;20)}
\def\zigzagCaa{\psline(2;0)(1.75;1)(2;2)(1.75;3)(2;4)(1.75;5)(2;6)(1.75;7)(2;8)(1.75;9)(2;10)}
\def\zigzagCbb{\psline(2;0)(1.87;1)(2;2)(1.87;3)(2;4)(1.87;5)(2;6)(1.87;7)(2;8)(1.87;9)(2;10)}
\rput{70}(0,0){\zigzagCa} \rput{20}(0,0){\zigzagCa}
\rput{-30}(0,0){\zigzagCa} \rput{120}(0,0){\zigzagCa}
\rput{45}(0,0){\zigzagCb} \rput{-5}(0,0){\zigzagCb}
\rput{95}(0,0){\zigzagCb}
\def\clocheCa{\psbezier(2;0)(2;0)(0.9;-15)(2;-30)}\def\clocheCb{\psbezier(2;0)(2;0)(0.3;-15)(2;-30)}
\rput{70}(0,0){\clocheCa} \rput{20}(0,0){\clocheCa}
\rput{120}(0,0){\clocheCa}
\rput{140}(0,0){\psline(2;0)(1.2;8)}\rput{-30}(0,0){\psline(2;0)(1.4;-10)}
\rput{45}(0,0){\clocheCb}\rput{95}(0,0){\clocheCb}
\rput{115}(0,0){\psline(2;0)(1.4;10)}\rput{-5}(0,0){\psline(2;0)(1.4;-10)}}
\end{pspicture}
\hspace{2cm}
\begin{pspicture}(-1.7,-1.7)(1.7,-1.7)
\scalebox{0.85}{ \psdots(0,2)
 \uput[ur](0,2){\scriptsize 1}
\SpecialCoor \pscircle(0,0){2}
\psarc[linestyle=dotted](0,0){1.5}{160}{-50}
\def\zigzagCa{\psline(2;0)(1.75;2)(2;4)(1.75;6)(2;8)(1.75;10)(2;12)(1.75;14)(2;16)(1.75;18)(2;20)}
\def\zigzagCb{\psline(2;0)(1.87;2)(2;4)(1.87;6)(2;8)(1.87;10)(2;12)(1.87;14)(2;16)(1.87;18)(2;20)}
\def\zigzagCaa{\psline(2;0)(1.75;1)(2;2)(1.75;3)(2;4)(1.75;5)(2;6)(1.75;7)(2;8)(1.75;9)(2;10)}
\def\zigzagCbb{\psline(2;0)(1.87;1)(2;2)(1.87;3)(2;4)(1.87;5)(2;6)(1.87;7)(2;8)(1.87;9)(2;10)}
\rput{70}(0,0){\zigzagCa} \rput{20}(0,0){\zigzagCa}
\rput{-30}(0,0){\zigzagCa} \rput{95}(0,0){\zigzagCa}
\rput{45}(0,0){\zigzagCb} \rput{-5}(0,0){\zigzagCb}
\rput{120}(0,0){\zigzagCb}
\def\clocheCa{\psbezier(2;0)(2;0)(0.9;-15)(2;-30)}\def\clocheCb{\psbezier(2;0)(2;0)(0.3;-15)(2;-30)}
\rput{70}(0,0){\clocheCa} \rput{20}(0,0){\clocheCa} 
\rput{140}(0,0){\psline(2;0)(1.2;8)}\rput{-30}(0,0){\psline(2;0)(1.4;-10)}
\rput{45}(0,0){\clocheCb}
\rput{115}(0,0){\psline(2;0)(1.4;10)}\rput{-5}(0,0){\psline(2;0)(1.4;-10)}
\psbezier(2;90)(2;90)(1.5;92.5)(2;95)
\psbezier(2;65)(2;65)(0.5;92.5)(2;120)}
\end{pspicture}
\caption{Circular representation of a
2-partition}\label{fig:Circular2partition}
\end{center}
\end{figure}
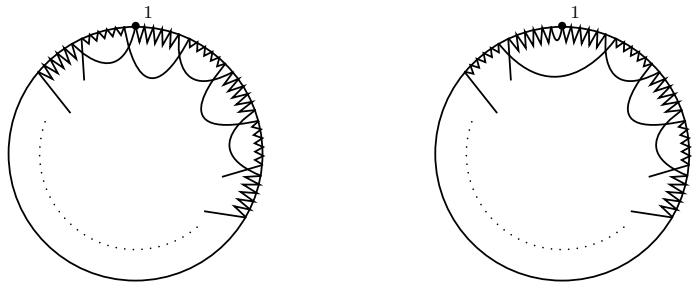

 A moment's thought (see Figure~\ref{fig:Circular2partition}) will convince
the reader that:
\begin{itemize}
\item if $k=1$, $\croc(\pi)=0$,
\item if $k\geq3$, $\croc(\pi)=2 k$.
\end{itemize}

This implies that, for $k\geq3$, we have
\begin{align*}
\sum_{\pi \in\T_k}
\croc(\pi)\,x^{\|\pi\|}=2k\,\left(\frac{x}{1-x}\right)^{2k}\,\frac{1}{1-x}
=2k\,\frac{x^{2k}}{(1-x)^{2k+1}},
\end{align*}
 since the generating function 
$G_k(x)=\sum_{\pi\in \T_k}x^{\|\pi\|}$ of 2-partitions in $\T_k$ is given by
$G_k(x)=\left(\frac{x}{1-x}\right)^{2k}\left(\frac{1}{1-x}\right)$.
Using the formal identity $\sum_{k\geq0}k\,y^{k}=\frac{y}{(1-y)^2}$,
it is then easy to obtain the identity
\begin{align}
\sum_{k\geq3}\sum_{\pi \in\T_k} \croc(\pi)\,x^{\|\pi\|}
&=\sum_{k\geq3}2k\,\frac{x^{2k}}{(1-x)^{2k+1}}
=\frac{2x^6(3-6x+x^2)}{(1-2x)^2(1-x)^5}.\label{eq:GfCircularTk}
\end{align}

\begin{figure}[h!]
\begin{center}
\small{
\begin{tabular}{|l||c|c|c|}
\cline{2-4}
\multicolumn{1}{c|}{}& & &\\[-0.3cm]
\multicolumn{1}{c|}{}& sketch of the circular & value & generating\\
\multicolumn{1}{c|}{}& representation & of $\croc$ & function\\
\multicolumn{1}{c|}{}& & &\\[-0.3cm]
\hline
& & &\\
& &  & \\
\quad case (a)& &&\\
$\ell_1=\ell_2=1$, $\ell_3=0$& & 1&$x^4$\\
 $m_1=m_2=1$ & &  & \\
& \begin{pspicture}(-1,-1.5)(1,-0.5)
 \scalebox{0.85}{\casA}
 \end{pspicture}& &\\[-0.3cm]
\hline
& & &\\
\quad case (b)& &&\\
$\ell_1=\ell_2=1$, $\ell_3=0$& &2& $x^2\left( \left(\frac{x} {1-x}\right)^{2}-x^2\right)$\\
  $m_1+m_2\geq3$ & &  & \\
&  \begin{pspicture}(-1,-1.5)(1,-0.5)
\scalebox{0.85}{\casB}
\end{pspicture} & &\\[-0.3cm]
\hline
& & &\\
\quad case (c)& &&\\
$\ell_1+\ell_2+\ell_3\geq3$& & 2& $x^2\left(\frac{x} {1-x}\right)^{2}\frac{1} {1-x}$\\
$m_1=m_2=1$& & &\\
&  \begin{pspicture}(-1,-1.5)(1,-0.5)
\scalebox{0.85}{\casCa}
\end{pspicture}
\begin{pspicture}(-1,-1.5)(1,-0.5)
\scalebox{0.85}{\casCb}
\end{pspicture}& & \\[-0.3cm]
\hline
& & &\\
\quad case (d)& &&\\
 $\ell_1+\ell_2+\ell_3\geq3$&&4&$\left(\frac{x} {1-x}\right)^{2}\frac{1} {1-x} \left( \left(\frac{x} {1-x}\right)^{2}-x^2\right)$\\
  $m_1+m_2\geq3$& & & \\
&  \begin{pspicture}(-1,-1.5)(1,-0.5)
\scalebox{0.85}{\casDa}
\end{pspicture}
\begin{pspicture}(-1,-1.5)(1,-0.5)
 \scalebox{0.85}{\casDb}
\end{pspicture} &  & \\
  \hline
\end{tabular}}
\end{center}
\caption{Circular crossings in a 2-partition in
$\T_2$.}\label{fig:crossingsInT2}
\end{figure}
 For the case $k=2$, we have to be more careful. By definition,
a partition in $\T_2$ can be written in the form
$\pi=1^{\ell_1}2^{m_1}1^{\ell_2}2^{m_2}1^{\ell_{3}}$, with
$\ell_1,\ell_2,m_1,m_2\geq1$ and $\ell_{3}\geq0$.   We distinguish
four cases. Instead of long-winded explanations, we prefer to present the results
in the table given in Figure~\ref{fig:crossingsInT2}. The first column
lists the different cases considered, while the second column indicates
 the corresponding circular representations. In each case,
the statistic $\croc$ assumes a constant value: the third column provides
these constants. Finally, the (ordinary) generating function of
2-partitions in each case is given in the fourth column.
 For instance, the first and third row read as follows:

$\bullet$ \emph{case (a)}:  there is one and only one  partition
$\pi=1\,3/2\,4$, $\croc(\pi)=1$ and
 the generating function for this case is $x^{4}$.

$\bullet$ \emph{case (c)}: the parameter $\croc$ is equal
to $2$  and the generating function of this class is
$x^2\left(\frac{x} {1-x}\right)^{2}\frac{1} {1-x}$.
 The other rows have to be read similarly.\\
Let $F_1$, $F_2$, $F_3$ and $F_4$ be the generating functions of 2-partitions in case $(a)$, $(b)$, $(c)$ and $(d)$,
respectively. It follows from a reading of the entire table that
\begin{align*}
\sum_{\pi \in\T_2} \croc(\pi)\,x^{\|\pi\|} = 1\,F_1(x)+2\,F_2(x)+
2\, F_3(x)+ 4\, F_4(x) =x^4+\frac{2x^5(5-4x+x^2)}{(1-x)^5}.
\end{align*}
Combining  \eqref{eq:GfCircularTk} and the latter equality, we
obtain
\begin{align*}
\sum_{\pi \in\Pi^2} \croc(\pi)\,x^{\|\pi\|}&=\sum_{k\geq1}\sum_{\pi
\in\T_k} \croc(\pi)\,x^{\|\pi\|}=\sum_{\pi \in\T_2} \croc(\pi)\,x^{\|\pi\|}+\sum_{k\geq3}\sum_{\pi \in\T_k} \croc(\pi)\,x^{\|\pi\|}\\
&=x^4+\frac{2x^5(5-4x+x^2)}{(1-x)^5}+\frac{2x^6(3-6x+x^2)}{(1-2x)^2(1-x)^5}\\
&=x^4+\frac{2x^5(5-11x+4x^2)}{(1-x)^3(1-2x)^2}.
\end{align*}

A routine coefficient extraction based on the partial fraction
decomposition
\begin{align*}
x^4+\frac{2x^5(5-11x+4x^2)}{(1-x)^3(1-2x)^2}&=-\frac{5 x}{2}-2
x^2+x^4-\frac{4}{(1-x)^3}+\frac{10}{(1-x)^2}-\frac{6}{1-x}\\
&\qquad +\frac{1}{4 (1-2 x)^2}-\frac{1}{4 (1-2 x)}
\end{align*}
 gives
$$
v_n=\sum_{\pi\in\Pi_n^2}\croc(\pi)=\begin{cases} 
                                         0, & \hbox{if $0\leq n\leq3$;} \\
                                         1, & \hbox{if $n=4$;} \\
                                         n\,2^{n-2}+4n -2n^2, & \hbox{if $n\geq5$.}
                                     \end{cases}
$$
A straightforward computation leads to
\begin{align}
 V_2^{\croc}(x)&=
\frac{1}{2}\,xe^{2x}+2x(-x+1)e^{x}+\frac{x^4}{24}-x^2-\frac{5}{2}\,x.\label{eq:gfV_CrossingCircular}
\end{align}
Since $V_2^{\croc}(x)$  can be written in the
form~\eqref{eq:canform_Vstat}, 
we can continue the process
described in Section~\ref{sect:process}.\\

\textit{Step 2.} Using~\eqref{eq:expgf_der12_Bell},
\eqref{eq:expgf_der1_Stirling} and \eqref{eq:expgf_der2_Stirling},
we can write
\begin{align*}
V_2^{\croc}(x)B(x)
&=\frac{1}{24} \left(x^4-24 x^2-60x\right)B(x)+ \frac{1}{2}(-4 x^2+3x)B^{\prime}(x)+\frac{x}{2} B^{\prime\prime}(x)\\
V_2^{\stat}(x)S_{k-2}(x)&= \frac{1}{24}
\left(x^4-72x^2\right)S_{k-2}(x)+ (k-1) \left(-2 x^2+3x\right)
S_{k-1}(x)\\
&\hspace{1cm}+\frac{1}{2} k(k-1) xS_{k}(x).
\end{align*}

\textit{Step 3.} After routine coefficient extraction based on
\eqref{eq:extractionCoeff}, we obtain
\begin{align}
 \left[\frac{x^{n}}{n!}\right]V_2^{\croc}(x)B(x)
&=\frac{n}{2}B_{n+1}+\frac{3n}{2}B_{n}
-2n\left(n+\frac{1}{4}\right)B_{n-1} -(n)_2B_{n-2}+\frac{1}{24}(n)_4
B_{n-4}\label{eq:V*B_croc}\\
 \left[\frac{x^{n}}{n!}\right]V_2^{\croc}(x)S_{k-2}(x)
& = \frac{1}{24}\left((n)_4 S_{n-4,k-2}-72(n)_2 S_{n-2,k-2}\right)\nonumber\\
 &\hspace{1cm}+ (k-1) \left(-2 (n)_2 S_{n-2,k-1}+3n S_{n-1,k-1}\right)\nonumber\\
&\hspace{1cm}+\frac{1}{2} k(k-1) nS_{n-1,k}.\nonumber
\end{align}

It is possible to simplify the expression obtained
 for $\left[\frac{x^{n}}{n!}\right]V_2^{\croc}(x)S_{k-2}(x)$ by replacing each occurrence
of the left hand sides of identities~\eqref{eq:StirlingSimplification4} and
\begin{align}
S_{n-2,k-1}&=\frac{1}{k-1}\left(S_{n-1,k-1}-S_{n-2,k-2}\right)
\end{align}
by the corresponding right hand sides.
This gives
\begin{align}
 \left[\frac{x^{n}}{n!}\right]V_2^{\croc}(x)S_{k-2}(x)
&=\frac{n(k-1)}{2}S_{n,k}+\frac{n}{2}(-4n-1+5k)S_{n-1,k-1}\label{eq:V*S_croc}\\
&\qquad -(n)_2 S_{n-2,k-2}+\frac{(n)_4}{24}S_{n-4,k-2}.\nonumber
 \end{align}

\textit{Steps 4 and 5.} Division of expression~\eqref{eq:V*B_croc}
by $B_n$ gives  the exact value~\eqref{eq:ExactMean_croc}
of $\mu_n$, while its asymptotic
approximation~\eqref{eq:AsymptoticMean_croc} is easily obtained
from Lemma~\ref{lem:Asymptotic_QuotientBell}.

 Similarly, division of expression~\eqref{eq:V*S_croc}
by $S_{n,k}$ gives the exact
value~\eqref{eq:ExactMeanBlock_croc} of $\mu_{n,k}$, while, after a
routine computation based on
Lemma~\ref{lem:Asymptotic_QuotientStirling}, we get the asymptotic
approximation
\begin{align*}
\mu_{n,k}&=\frac{1}{2}n(k-1)+ O \left(
n^2\left(1-\frac{1}{k}\right)^n \right),\quad n\rightarrow\infty,
\end{align*}
which is a refinement of~\eqref{eq:AsymptoticMeanBlock_croc}. This
concludes the proof of Theorem~\ref{thm:Mean_CrossingCircular}.

\subsection{Some Remarks}
 While the combinatorial parameter $\crol$ has received a
considerable interest, it seems that the natural parameter~$\croc$
was ignored. It is interesting to ask whether this parameter admits ``closed''
formulas similar to the formulas for $\crol$ presented in the previous section.
 Some results on the asymptotic distribution of the parameter~$\croc$
will be given in~\cite{KaFutur}.


\section{Number of overlappings}

\subsection{Proof of Theorem~\ref{thm:exactmean_ovl}}

\emph{Computation of $V_2^{\ov}(x)$.} By definition,  $V_2^{\ov}(x)$
is the exponential generating function of the sequence
$\left(v_n\right)_{n\ge0}$ with
$v_n:=v^{\ov}_{n,2}=\sum_{\pi\in\Pi_n^2}\ov(\pi)$.

 Since $\ov$ assumes only the values 0 or 1 on $\Pi_n^{2}$, the term $v_n$ is just
the number of $2$-partitions of $[n]$ that overlap, i.e., the number
of  2-partitions $\pi=B_1/B_2$ of $[n]$ such that 
$1=\min B_1< \min B_2 <\max B_1< \max B_2=n$. It is easy
to see that there are exactly $2^{n-2}$ partitions $\pi=B_1/B_2$ of
$[n]$ such that $\min B_1=1$ and $\max B_2=n$, and among them only
$n-1$ do not overlap (these are the partitions $\pi=1^{k}2^{n-k}$
with $1\leq k \leq n-1$). Consequently
$$
v_n=\sum_{\pi\in\Pi_n^2}\ov(\pi)=\begin{cases} 
                                         2^{n-2}-n+1, & \hbox{if $n\geq2$;} \\
                                         0, & \hbox{if $n\leq 1$.}
                                    \end{cases}
$$

 A straightforward computation leads to
\begin{align}\label{eq:gfVov}
 V_2^{\ov}(x)=\sum_{n\geq0}v_n\,\frac{x^n}{n!}&=\frac{1}{4}e^{2x}+(-x+1)\,e^x-\frac{x}{2}-\frac{5}{4}.
\end{align}
Since $V_2^{\ov}(x)$  can be written in the
form~\eqref{eq:canform_Vstat}, 
we can continue the process
described in Section~\ref{sect:process}.\\

\textit{Step 2.} Using~\eqref{eq:expgf_der12_Bell},
\eqref{eq:expgf_der1_Stirling} and \eqref{eq:expgf_der2_Stirling},
we can write
\begin{align*}
V_2^{\ov}(x)B(x)
&= \left(-\frac{x}{2}-\frac{5}{4}\right)B(x)+\left(-x+\frac{3}{4}\right)B^{\prime}(x)+\frac{1}{4} B^{\prime\prime}(x)\\
V_2^{\ov}(x)S_{k-2}(x)&= -\frac{3 x}{2} S_{k-2}(x)-
(k-1)\left(x-\frac{3}{2}\right)S_{k-1}(x)+ \frac{1}{4}  k(k-1)
S_{k}(x).
\end{align*}

\textit{Step 3.} After routine coefficient extraction based on
\eqref{eq:extractionCoeff}, we obtain
\begin{align}
 \left[\frac{x^{n}}{n!}\right]V_2^{\ov}(x)B(x) &=\frac{1}{4}
B_{n+2}+\frac{3}{4} B_{n+1}-\left(n+\frac{5}{4}\right)
B_{n}-\frac{n}{2}B_{n-1}\label{eq:V*B_ov}\\
 \left[\frac{x^{n}}{n!}\right]V_2^{\ov}(x)S_{k-2}(x) &=
-\frac{3n}{2} S_{n-1,k-2}-\frac{1}{2} (k-1)\left(2
nS_{n-1,k-1}-3S_{n,k-1}\right) \nonumber\\
&\hspace{1cm}+ \frac{1}{4}  k(k-1) S_{n,k}.\nonumber
\end{align}
It is possible to simplify the expression obtained
 for $\left[\frac{x^{n}}{n!}\right]V_2^{\ov}(x)S_{k-2}(x)$ by replacing each occurrence
of the left hand side  of identity~\eqref{eq:StirlingSimplification2}
by the corresponding right hand side.
This gives
\begin{align}
 \left[\frac{x^{n}}{n!}\right]V_2^{\ov}(x)S_{k-2}(x)
&=\frac{1}{4}k(k-1)S_{n,k}-\frac{3}{2}(n+1-k)S_{n,k-1} \label{eq:V*S_ov}\\
&\qquad +\frac{1}{2}n(k-1)S_{n-1,k-1}.\nonumber
\end{align}

\textit{Steps 4 and 5.} Division of expression~\eqref{eq:V*B_ov} by
$B_n$ gives  the exact value~\eqref{eq:ExactMean_ov} of
$\mu_n$, while its asymptotic
approximation~\eqref{eq:AsymptoticMean_ov} is easily obtained from
Lemma~\ref{lem:Asymptotic_QuotientBell}.

 Similarly, division of expression~\eqref{eq:V*S_ov}
by $S_{n,k}$ gives the exact
value~\eqref{eq:ExactMeanBlock_ov} of $\mu_{n,k}$, while, after a
routine computation based on
Lemma~\ref{lem:Asymptotic_QuotientStirling}, we obtain the asymptotic
approximation
\begin{align*}
\mu_{n,k}&=\frac{1}{2}{k\choose2}+ O \left(
n\left(1-\frac{1}{k}\right)^n \right),\quad n\rightarrow\infty,
\end{align*}
which is a refinement of~\eqref{eq:AsymptoticMeanBlock_ov}. This
concludes the proof of Theorem~\ref{thm:exactmean_ovl}.

\subsection{Some remarks}
Although it seems that the parameter number of overlappings was never considered before,
it is naturally suggested by the article~\cite{FlS}, in which Flajolet and Schott
considered non-overlappings partitions. In particular, they obtained the generating function of
non-overlapping partitions in the form of a continued fraction expansion.
 It is worth noting that a continued fraction expansion for 
the generating function of set partitions with respect to the number
of overlappings follows painlessly from earlier combinatorial
investigations (see e.g.~\cite{Fl,FlS}). 

We also want to point out that, as $n\rightarrow\infty$, $\mu_{n,k}$
converges to ${k\choose 2}$. That this is also the average number of inversions in
a random permutation of size~$k$ is far from being a coincidence, as we will
explain in future work.


\section{Number of Embracings}

\subsection{Proof of Theorem~\ref{thm:exactmean_emb}}

\emph{Computation of $V_2^{\emb}(x)$.} By definition,
$V_2^{\emb}(x)$ is the exponential generating function of the
sequence $\left(v_n\right)_{n\ge0}$ with
$v_n:=v^{\emb}_{n,2}=\sum_{\pi\in\Pi_n^2}\emb(\pi)$.

 Since $\emb$ assumes only the values 0 or 1 on $\Pi_n^{2}$, the term $v_n$ is just the number
of 2-partitions of~$[n]$ which are embracing, i.e., the number of
2-partitions $\pi=B_1/B_2$ of~$[n]$ such that $1=\min B_1< \min B_2
\leq \max B_2< \max B_1=n$, or, equivalently, such that $\min B_1=1$
and $\max B_1=n$. Obviously, for $n\geq 2$, there are exactly
$2^{n-2}-1$ such 2-partitions. Consequently
$$
v_n=\sum_{\pi\in\Pi_n^2}\emb(\pi)=\begin{cases} 
                                         2^{n-2}-1, & \hbox{if $n\geq2$;} \\
                                         0, & \hbox{if $n\leq 1$.}
                                     \end{cases}
$$
 A straightforward computation leads to
\begin{align}\label{eq:gfVemb}
 V_2^{\emb}(x)=\sum_{n\geq0}v_n\,\frac{x^n}{n!}&=\frac{1}{4}e^{2x}-\,e^x+\frac{x}{2}+\frac{3}{4}.
\end{align}
Since $V_2^{\emb}(x)$  can be written in the
form~\eqref{eq:canform_Vstat}, 
we can continue the process
described in Section~\ref{sect:process}.\\

\textit{Step 2.} Using~\eqref{eq:expgf_der12_Bell},
\eqref{eq:expgf_der1_Stirling} and \eqref{eq:expgf_der2_Stirling},
we can write
\begin{align*}
V_2^{\emb}(x)B(x)
&= \left(\frac{x}{2}+\frac{3}{4}\right)B(x)-\frac{5}{4} B^{\prime}(x)+\frac{1}{4} B^{\prime\prime}(x)\\
V_2^{\emb}(x)S_{k-2}(x)&= \frac{x}{2} S_{k-2}(x)-\frac{1}{2}
(k-1)S_{k-1}(x)+ \frac{1}{4}  k(k-1) S_{k}(x).
\end{align*}

\textit{Step 3.} After routine coefficient extraction based on
\eqref{eq:extractionCoeff}, we obtain
\begin{align}
 \left[\frac{x^{n}}{n!}\right]V_2^{\emb}(x)B(x) &=\frac{1}{4}
B_{n+2}-\frac{5}{4}
B_{n+1}+\frac{3}{4}B_{n}+\frac{n}{2}B_{n-1}\label{eq:V*B_emb}\\
 \left[\frac{x^{n}}{n!}\right]V_2^{\emb}(x)S_{k-2}(x) &=
\frac{n}{2} S_{n-1,k-2}-\frac{1}{2} (k-1)S_{n,k-1}+
\frac{1}{4}k(k-1)S_{n,k}.\label{eq:V*S_emb}
\end{align}

\textit{Steps 4 and 5.} Division of expression~\eqref{eq:V*B_emb}
by $B_n$ gives the exact value~\eqref{eq:ExactMean_emb}
of $\mu_n$, while its asymptotic
approximation~\eqref{eq:AsymptoticMean_emb} is easily obtained from
Lemma~\ref{lem:Asymptotic_QuotientBell}.

 Similarly, division of expression~\eqref{eq:V*S_emb}
by $S_{n,k}$ gives the exact
value~\eqref{eq:ExactMeanBlock_emb} of $\mu_{n,k}$, while, after a
routine computation based on
Lemma~\ref{lem:Asymptotic_QuotientStirling}, we get the asymptotic
approximation
\begin{align*}
\mu_{n,k}&=\frac{1}{2}{k\choose2}+ O \left(
\left(1-\frac{1}{k}\right)^n \right),\quad n\rightarrow\infty,
\end{align*}
which is a refinement of~\eqref{eq:AsymptoticMeanBlock_emb}. This
concludes the proof of Theorem~\ref{thm:exactmean_emb}.

\subsection{Some remarks}
It is worth noting that, like for overlappings, 
it is painless to derive a continued fraction form for the generating function of set
partitions with respect
to the number of embracings by making use of the combinatorial theory
of continued fractions~\cite{Fl}.

 One can also consider a variation of embracings. Say that two
 sets~$B$ and~$B'$ \emph{strongly embrace each other}
if, with the implied order structure, $\min(B)<\min(B')<\max(B')< \max(B)$. Then one can show that
the statistic `number of strong embracings' has the same distribution  (hence, the same average value) on each $\Pi_n^k$
as the number of overlappings. This can be proved using the combinatorial theory
of continued fractions or by a direct combinatorial argument.


\section{Number of occurrences of a 2-Pattern}

\subsection{Proof of Theorem~\ref{thm:exactmean_occ} and Theorem~\ref{thm:exactmeanblock_occ} }

\emph{Computation of $V_2^{\occs}(x)$.} 
Let $\sigma=\sigma_1\sigma_2\ldots\sigma_r$ be a $2$-pattern of length $r$.
By definition, $V_2^{\occs}(x)$ is the exponential generating function of the
sequence $\left(v_n\right)_{n\ge0}$ with
${v_n:=v^{\occs}_{n,2}=\sum_{\pi\in\Pi_n^2}\occs(\pi)}$.

We can interpret $v_n$ as the number of pairs $(\pi,t)$,
where~$\pi$ is a 2-partition of $[n]$ and~$t=(i_1,i_2,\ldots,i_r)$
is an occurrence of~$\sigma$ in~$\pi$. Such a pair will be called an
\emph{underlined} 2-partition.
 An underlined 2-partition $(\pi,t)$
with $t=(i_1,i_2,\ldots,i_r)$ can be identified with the restricted
growth function $w(\pi)$ in which the letters
$w_{i_1}$,$w_{i_2}$,\ldots,$w_{i_r}$ are colored. For instance, if
$\sigma=2\,1\,2$, the 3-tuple $t=(2,4,5)$ is an occurrence of
$\sigma$ in $\pi=1\,4/2\,3\,5$ and
$(\pi,t)\equiv 1\,\mathbf{2}\,2\,\mathbf{1}\,\mathbf{2}$.

 Using this correspondence, it is not hard
to see that an underlined 2-partition can be decomposed uniquely as
\begin{align*}
\pi=u_1 \, \mathbf{\sigma_{1}}\,u_2 \, \mathbf{\sigma_{2}}\cdots
\,u_r \, \mathbf{\sigma_{r}}\,u_{r+1}\quad\text{if $\sigma_1=1$},
\end{align*}
with $u_1\in\{\epsilon\}\cup 1\{1,2\}^*$ and $u_i\in\{1,2\}^*$ for
$i=2,\ldots,r+1$;
\begin{align*}
\pi=1\,v_1 \, \mathbf{\sigma_{1}}\,v_2 \, \mathbf{\sigma_{2}}\cdots
\,v_r \, \mathbf{\sigma_{r}}\,v_{r+1}\quad\text{if $\sigma_1=2$},
\end{align*}
with $v_i\in\{1,2\}^*$ for $i=1,\ldots,r+1$. In particular, this
proves that $v_n$ (and thus $V_2^{\occs}(x)$) depends only on the
first letter of $\sigma$. For $i\in\{1,2\}$, let
$v_n^{(\sigma_1=i)}$ equal the value of $v_n$ when $\sigma_1=i$. It
follows from the above decomposition of underlined 2-partitions 
combined with an elementary counting that for $n\geq r$,
$$
v_n=\sum_{\pi\in\Pi_n^2}\occs(\pi)=\begin{cases} 
                                        v_n^{(\sigma_1=1)}={n-1 \choose r-1}\,2^{n-r}+{n-1 \choose r}\,2^{n-r-1}, & \hbox{if $\sigma_1=1$;} \\
                                        v_n^{(\sigma_1=2)}= {n-1 \choose r}\,2^{n-r-1}, & \hbox{if $\sigma_1=2$.}
                                    \end{cases}
$$
To see the above result in the case $\sigma_1=2$, for instance, just
observe that there are ${n-1 \choose r}$ choices for the positions of the
underlined elements
$\mathbf{\sigma_1}$,$\mathbf{\sigma_2}$,\ldots,$\mathbf{\sigma_r}$
and then $2^{n-r-1}$ choices for the word $v_1 \, v_2 \cdots
\,v_{r+1}$. The case $\sigma_1=1$ can be treated in a similar way.

For $i\in\{1,2\}$, let $V_2^{(\sigma_1=i)}(x)$ equal
$V_2^{\occs}(x)$ in the case where $\sigma_1=i$, i.e.,
$V_2^{(\sigma_1=i)}(x)=\sum_{n\geq
0}v_n^{(\sigma_1=i)}\,\frac{x^n}{n!}$. For later simplifications, it
is important to note that it suffices to determine only one of the
$V_2^{(\sigma_1=i)}(x)$'s to get the other. Indeed, an elementary
computation yields
\begin{align*}
 v_n^{(\sigma_1=1)}+v_n^{(\sigma_1=2)}
=\left({n-1 \choose r-1}\,2^{n-r}+{n-1 \choose r}\right)\,2^{n-r-1}+{n-1 \choose r}\,2^{n-r-1}
={n \choose r}\,2^{n-r},
\end{align*}
from which we deduce easily (using \eqref{eq:DeriveesFormelles}, for instance) that
\begin{align}\label{eq:V1+V2}
 V_2^{(\sigma_1=1)}(x)+V_2^{(\sigma_1=2)}(x)
= \sum_{n\geq 0}\left( v_n^{(\sigma_1=1)}+v_n^{(\sigma_1=2)}\right)\,\frac{x^n}{n!}
=\sum_{n\geq 0} {n \choose r}\,2^{n-r}\,\frac{x^n}{n!}
=\frac{x^r}{r!}\,e^{2x}.
\end{align}

We choose to deal with $V_2^{(\sigma_1=2)}(x)$. In order to obtain a convenient expression
 for  $V_2^{(\sigma_1=2)}(x)$, we
need some additional materials. Let $k$ be a nonnegative integer and
${(n)_k:=n(n-1)(n-2)\cdots(n-k+1)}$, as before.
Given a formal power series $A(x):=\sum_{n\geq0}a_n\frac{x^n}{n!}$,
it is easy to establish the formal identity
\begin{align}\label{eq:DeriveesFormelles}
\sum_{n\geq0}(n)_k a_n\frac{x^n}{n!}=x^kA^{(k)}(x).
\end{align}
What about the power series $\sum_{n\geq0}(n-1)_k\,
a_n\frac{x^n}{n!}$? Using the Chu--Vandermonde convolution, we have
\begin{align*}
(n-1)_k=\sum_{j=0}^k {k\choose j}(n)_j (-1)_{k-j}=\sum_{j=0}^k
(-1)^{k-j}\frac{k!}{j!}(n)_j,
\end{align*}
from which we deduce that
\begin{align}
\sum_{n\geq0}(n-1)_k \,a_n\frac{x^n}{n!}=\sum_{j=0}^k
(-1)^{k-j}\frac{k!}{j!}x^jA^{(j)}(x).
\end{align}
 In particular, for $a_n=2^{n}$, we obtain
\begin{align}
\sum_{n\geq0}(n-1)_k 2^n\frac{x^n}{n!}=\sum_{j=0}^k
(-1)^{k-j}\frac{k!}{j!}x^j2^j e^{2x}=\left(\sum_{j=0}^k
(-1)^{k-j}\frac{k!}{j!}x^j2^j \right)
e^{2x}.\label{eq:gen_ShiftedLowFactorial}
\end{align}

  We can now express the power series $V_2^{(\sigma_1=2)}(x)$ (and   $V_2^{(\sigma_1=1)}(x)$) in a convenient form.
We have seen earlier that $v_n^{(\sigma_1=2)}={n-1 \choose
r}\,2^{n-r-1}$ for $n\geq r$ ($v_n^{(\sigma_1=2)}=0$, otherwise).
After elementary manipulations, one gets
\begin{align}
V_2^{(\sigma_1=2)}(x)&=\sum_{n\geq r}{n-1 \choose
r}\,2^{n-r-1}\frac{x^n}{n!}
=\frac{(-1)^{r+1}}{2^{r+1}}+\frac{1}{r!\,2^{r+1}}\sum_{n\geq
0}(n-1)_r\,2^n\frac{x^n}{n!}.\label{eq:Vocc2_intermed}
\end{align}
For any integer $r\geq 0$, set
\begin{align}\label{eq:defPr(x)}
P_r(x):=\sum_{j=0}^r\frac{(-1)^{r-j}}{j!2^{r+1-j}}x^j.
\end{align}
Combining~\eqref{eq:Vocc2_intermed}
and~\eqref{eq:gen_ShiftedLowFactorial}, we arrive at
\begin{align}\label{eq:Vocc2}
V_2^{(\sigma_1=2)}(x)&=\frac{(-1)^{r+1}}{2^{r+1}}+P_r(x)e^{2x}.
\end{align}
Combining the latter identity and~\eqref{eq:V1+V2}, we  get
\begin{align}\label{eq:Vocc1}
V_2^{(\sigma_1=1)}(x)&=\frac{(-1)^{r}}{2^{r+1}}+\left(-P_r(x)+\frac{x^r}{r!}\right)\,e^{2x}.
\end{align}

\vspace{0.3cm}

\emph{The connecting relations.} We now establish
relations~\eqref{eq:ConnectingRelation_occs}
and~\eqref{eq:ConnectingRelationBlock_occs}. Instead
of~\eqref{eq:expgf_der2_Stirling}, we will use the identity
\begin{align}
e^{2x}S_{k-2}(x)&=S^{\prime\prime}_{k}(x)-S^{\prime}_{k}(x),\label{eq:expgf_der3_Stirling}
\end{align}
which can be obtained from~\eqref{eq:Stirling_VerticalGF} by a
straightforward computation.  Combining~\eqref{eq:V1+V2}
and \eqref{eq:expgf_der3_Stirling}, we can write
\begin{align*}
\left(V_2^{(\sigma_1=1)}(x)+V_2^{(\sigma_1=2)}(x)\right)\,S_{k-2}(x)
&=\frac{x^r}{r!}\,\left(S^{\prime\prime}_{k}(x)-S^{\prime}_{k}(x)\right),
\end{align*}
and, after a routine coefficient extraction based on
\eqref{eq:extractionCoeff},  we obtain {\small
\begin{align}\label{eq:ConnectingRelation_occs_intermed}
 \left[\frac{x^{n}}{n!}\right]\left(V_2^{(\sigma_1=1)}(x)+V_2^{(\sigma_1=2)}(x)\right)\,S_{k-2}(x)
&={n \choose r} \left(S_{n+2-r,k}-S_{n+1-r,k}\right).
\end{align}}
Summing the both sides of this identity over all integers $k\geq 0$, we obtain
\begin{align}\label{eq:ConnectingRelationBlock_occs_intermed}
 \left[\frac{x^{n}}{n!}\right]\left(V_2^{(\sigma_1=1)}(x)+V_2^{(\sigma_1=2)}(x)\right)\,B(x)
&={n \choose r} \left(B_{n+2-r}-B_{n+1-r}\right).
\end{align}

It suffices now to divide
expression~\eqref{eq:ConnectingRelation_occs_intermed} by~$B_n$ and
expression~\eqref{eq:ConnectingRelationBlock_occs_intermed} by
$S_{n,k}$, respectively, to obtain the
relations~\eqref{eq:ConnectingRelation_occs}
and~\eqref{eq:ConnectingRelationBlock_occs}.\\

\emph{Exact average values in the case $\sigma_1=2$.}
Combining~\eqref{eq:Vocc2} and~\eqref{eq:expgf_der3_Stirling}, we
can write
\begin{align*}
 V_2^{(\sigma_1=2)}S_{k-2}(x)
 &=\frac{(-1)^{r+1}}{2^{r+1}}S_{k-2}(x)+P_r(x) S^{\prime\prime}_{k}(x)-P_r(x) S^{\prime}_{k}(x).
 \end{align*}
After a routine coefficient extraction based on
\eqref{eq:extractionCoeff}, we obtain 
\begin{align*}
\left[\frac{x^{n}}{n!}\right]V_2^{(\sigma_1=2)}(x)S_{k-2}(x)&=\sum_{j=0}^r\frac{(-1)^{r-j}}{j!2^{r+1-j}}(n)_j
\left(S_{n+2-j,k}-S_{n+1-j,k}\right)\\&\qquad +
\frac{(-1)^{r+1}}{2^{r+1}}S_{n,k-2}.
\end{align*}
A rearrangement of the terms of the sum on the right hand side gives
\begin{align}
\left[\frac{x^{n}}{n!}\right]V_2^{(\sigma_1=2)}(x)S_{k-2}(x)&=\frac{(-1)^{r+1}}{2^{r+1}}S_{n,k-2}+\sum_{j=0}^r p_j(n) S_{n+2-j,k}\label{eq:V*S_occ2}\\
 &\qquad -\frac{1}{2}{n\choose r} S_{n+1-r,k},\nonumber
\end{align}
where $p_j(n)$ is given by~\eqref{eq:defP_j}.

Division of expression~\eqref{eq:V*S_occ2} by $S_{n,k}$ gives the
exact value of $\mu_{n,k}^{(\sigma_1=2)}$ given in
Theorem~\ref{thm:exactmeanblock_occ}, while the exact average value
of $\mu_n^{(\sigma_1=2)}$ given in
Theorem~\ref{thm:exactmeanblock_occ} can be obtained easily from
the exact value of $\mu_{n,k}^{(\sigma_1=2)}$ by
using~\eqref{eq:Meanblock_Mean}.\\

\emph{Asymptotic approximations of the average values in the case
$\sigma_1=2$.} The exact average value of $\mu_{n}^{(\sigma_1=2)}$
we found is
\begin{align*}
\mu_n^{(\sigma_1=2)}&=\frac{(-1)^{r+1}}{2^{r+1}}+\sum_{j=0}^r p_j(n)\frac{B_{n+2-j}}{B_n}
-\frac{1}{2}{n\choose r}\frac{B_{n+1-r}}{B_n}.
\end{align*}
It follows easily from~Lemma~\ref{lem:Asymptotic_QuotientBell} and
basic representations of binomial coefficients that, as
$n\rightarrow\infty$, we have
\begin{align}
p_j(n)\frac{B_{n+2-j}}{B_n}&=\frac{(-1)^{r-j}}{2^{r+1-j}}\left({n\choose j}+\frac{1}{2}{n\choose
j-1}\right)\frac{B_{n+2-j}}{B_n} \nonumber\\
 &=\frac{(-1)^{r-j}}{2^{r+1-j}}
\frac{n^{j}}{j!}\left(1+O\left(\frac{1}{n}\right)\right)  \left(
\frac{n}{\log n}\right)^{2-j}
\left( 1+(2-j) \frac{\log\log n}{\log n} \left(1+o\left(1\right)\right)  \right)\nonumber \\
&=\frac{(-1)^{r-j}}{2^{r+1-j}} \frac{1}{j!} n^{2} \left(\log n\right)^{j-2} \left( 1-(j-2)
\frac{\log\log n}{\log n} +o\left(\frac{\log\log n}{\log n}\right) \right)\label{eq:asymptotic_occ_lem1}
\end{align}
for any integer $j\geq 0$. Similarly, as $n\rightarrow\infty$, one
can prove
\begin{align}
{n\choose r}\frac{B_{n+1-r}}{B_n} &=
 \frac{1}{r!} n \left(\log n\right)^{r-1} \left( 1-(r-1)
\frac{\log\log n}{\log n} \left(1+o\left(1\right)\right)  \right).\label{eq:asymptotic_occ_lem2}
\end{align}
 Combining~\eqref{eq:asymptotic_occ_lem1} and~\eqref{eq:asymptotic_occ_lem2}, we get
after basic comparisons
\begin{align*}
\mu_n^{(\sigma_1=2)}&=\frac{1}{2\,r!} n^{2} \left(\log
n\right)^{r-2} \left( 1-(r-2) \frac{\log\log n}{\log n}
+o\left(\frac{\log\log n}{\log n}\right) \right),
\end{align*}
which is exactly~\eqref{eq:AsymptoticMean_occ} in the case
$\sigma_1=2$.
\vspace{0.3cm}

The exact average value of $\mu_{n,k}^{(\sigma_1=2)}$ we found is
\begin{align*}
\mu_{n,k}^{(\sigma_1=2)}&=\sum_{j=0}^rp_j(n)\frac{S_{n+2-j,k}}{S_{n,k}}-\frac{1}{2}{n\choose r}
\frac{S_{n+1-r,k}}{S_{n,k}}+ \frac{(-1)^{r+1}}{2^{r+1}}\frac{S_{n,k-2}}{S_{n,k}}.
\end{align*}
It follows easily from~Lemma~\ref{lem:Asymptotic_QuotientStirling}
and basic representations of binomial coefficients that, as
$n\rightarrow\infty$, we have
\begin{align*}
\mu_{n,k}^{(\sigma_1=2)}&=\sum_{j=0}^r\frac{(-1)^{r-j}}{2^{r+1-j}}
\left({n\choose j}+\frac{1}{2}{n\choose j-1}\right)
 \left( k^{2-j} +O\left(\left(1-\frac{1}{k}\right)^n\right)  \right)\\
& \qquad -\frac{1}{2}{n\choose r} \left( k^{1-r}
+O\left(\left(1-\frac{1}{k}\right)^n\right)  \right)
 + \frac{(-1)^{r+1}}{2^{r+1}}O\left(\left(1-\frac{2}{k}\right)^n\right)\\
&=\sum_{j=0}^r\frac{(-1)^{r-j}}{2^{r+1-j}}\left({n\choose
j}+\frac{1}{2}{n\choose j-1}\right)k^{2-j}\;
 -\frac{1}{2}{n\choose r}k^{1-r} +O\left(n^{r}\left(1-\frac{1}{k}\right)^n\right),
\end{align*}
which is a refinement of the asymptotic
approximation~\eqref{eq:AsymptoticMeanBlock_occ2} given in
Theorem~\ref{thm:exactmeanblock_occ}.\\

\emph{The case $\sigma_1=1$.} The exact and asymptotic approximations
of the average value $\mu_{n}^{(\sigma_1=1)}$ ($\mu_{n,k}^{(\sigma_1=1)}$,
respectively) given in Theorem~\ref{thm:exactmean_occ}
(Theorem~\ref{thm:exactmeanblock_occ},
respectively) can be  easily
obtained from the exact and asymptotic approximations of
$\mu_{n}^{(\sigma_1=2)}$ ($\mu_{n,k}^{(\sigma_1=2)}$,
respectively)
obtained earlier combined with 
relation~\eqref{eq:ConnectingRelation_occs}
(respectively \eqref{eq:ConnectingRelationBlock_occs}). The details are
left to the reader.\\

\emph{Patterns of length $2$.} To end the proof of our theorems, we
need to establish assertion $(v)$ of
Theorem~\ref{thm:exactmean_occ}. In fact, we have already shown in
Section~\ref{section:inversion}
that~\eqref{eq:AsymptoticMean_occ_length2} is an asymptotic
approximation of $\mu_n^{(2\,1)}$ (which is actually equal to
$\mu_n^{\inv}$ as we remarked earlier). Almost the same asymptotic
computation shows that~\eqref{eq:AsymptoticMean_occ_length2} is an
asymptotic approximation of $\mu_n^{(1\,2)}$ whose exact value is
given in~\eqref{eq:mean12} (the details are left to the reader).

This concludes the proof of Theorem~\ref{thm:exactmean_occ} and
Theorem~\ref{thm:exactmeanblock_occ}.

\vspace{0.3cm}

\subsection{Remarks}\label{section:remarks_occ}

 The notion of pattern containment in set partitions we have considered 
in this section is due to Sagan~\cite{Sag1}. There is also an older and 
different notion of pattern containment in set partitions due to Klazar 
(see e.g.~\cite{Kl}). Actually, it can be shown that the statistic `number 
of occurrences of a 2-partition' (according to Klazar's definition) is also 
a Z-statistic. One can show by either applying the methodology
developed in this paper or by direct combinatorial arguments
that, with Klazar's definition of pattern containment, the average
number  of occurrences of any 2-partition $\pi$ of
$[r]$ in a random set partition of $[n]$ is given by
\begin{align*}
  {n \choose r} \frac{B_{n+2-r}-B_{n+1-r}}{B_{n}}.
\end{align*}

 Mansour, Shattuck and Yan~\cite{Mans} have considered a restricted case of
pattern containment in set partitions where it is required that the
occurrence has to be contiguous. They obtained some average value
results for  numbers of restricted occurrences of some patterns. Note
that their results can not be recovered with the methodology of this
paper since, in the contiguous case, numbers of occurrences are not
Z-statistics.


\section{Concluding remarks}

\subsection{Average values of Z-statistics of depth $r$}

 There is a natural way to generalize Z-statistics. Given an integer $r\geq 1$,
we will say that a set partition statistic $STAT$ is a Z-statistic of depth $r$ if for any $\pi\in\Pi$,
\begin{align}
 STAT(\pi)=\sum_{A_1,A_2,\ldots,A_r\in \pi} STAT(\st(A_1/A_2/\cdots/A_r)).\label{eq:defn_rsplittable_prop1}
 \end{align}
Examples of Z-statistics of depth $r$ which can be found in the literature are the numbers of $r$-crossings and $r$-nestings
of arcs~\cite{ChenSt}, the number of occurrences of an $r$-pattern
with Sagan's definition~(see Section~\ref{section:remarks_occ}) and the number of occurrences of an $r$-partition
with Klazar's definition~(see Section~\ref{section:remarks_occ}).

  Theorem~\ref{thm:GfunMean} can be generalized  in a straightforward manner to simplify the computation of
average values of Z-statistics of depth $r$.
 Given a Z-statistic of depth $r$  $\stat$, let $v^{\stat}_{n,r}:=\sum_{\pi\in\Pi_n^r}\stat(\pi)$, and
let $V_r^{\stat}(x)$ be the exponential generating function of the
sequence $\left(v^{\stat}_{n,r}\right)_{n\ge0}$, i.e.,
\begin{align}
 V_r^{\stat}(x)=\sum_{n\geq0} v^{\stat}_{n,r}\,\frac{x^{n}}{n!}=\sum_{n\geq0} \sum_{\pi\in\Pi_n^r}\stat(\pi)\,\frac{x^{n}}{n!}.
\end{align}

Then one can prove the following result.
\begin{thm}\label{thm:GfunMean_r-splittable}
Let $\stat$ be  a Z-statistic of depth $r$. Denote by  $\mu_n$ (respectively
$\mu_{n,k}$)  the average value of the statistic $\stat$ in a random
set partition of~$[n]$ (respectively of~$[n]$ into $k$ blocks). Then we
have:
\begin{align}
\mu_n&=\frac{1}{B_n}\,\left[\frac{x^n}{n!}\right]\,V_r^{\stat}(x)B(x)\label{eq:ThmGfunMean_r-split}\\
\mu_{n,k}&=\frac{1}{S_{n,k}}\,\left[\frac{x^n}{n!}\right]\,V_r^{\stat}(x)S_{k-2}(x).\label{eq:ThmGfunMeanBlock_r-split}
\end{align}
\end{thm}
Note that the point of the above result is that the computation
of average values of a Z-statistic of depth $r$ in a random set
partition is essentially equivalent to the computation of its
average value in a random $r$-partition.

 As open problems, we ask the following questions.
What are the (exact and asymptotic) average values of the numbers of
$r$-crossings, $r$-nestings, occurrences of an $r$-pattern in a
random set partition of $[n]$?

\subsection{Other exponential families}

 In this paper we concentrated on set partitions but we might consider any other exponential family.
An exponential family is, roughly speaking, a set of structures that are built out of
connected pieces  (see e.g.~\cite{Wilf} for a precise definition). One can easily establish a generalized version
of Theorem~\ref{thm:GfunMean} for exponential families to simplify computation of average values.
We do not think that it is necessary to state the generalization, but we give a typical illustration of application.

Consider a natural generalization of matchings, namely the
(exponential family of) set partitions all blocks of which are
of size $m$ ($m$ is a given integer $\geq 1$). 
For these partitions, 
we can consider the statistics `numbers of crossings in  linear  and circular representations'
and it is easy to see that these statistics are still Z-statistics on such partitions. Then one can prove
the following results.

$\bullet$ \emph{The average number $\mu^{(\ell)}_{(m^k)}$ of linear crossings
in a random set partition into $k$~blocks of size~$m$, $m\geq1$, is
given by}
\begin{align*}
\mu_{(m^k)}={k\choose 2}\left(m-2+\frac{2}{{2m\choose m}}\right).
\end{align*}

$\bullet$ \emph{The average number $\mu^{(c)}_{(m^k)}$ of circular crossings
in a random set partition into $k$~blocks of size~$m$, $m\geq 3$, is
given by}
\begin{align*}
\mu_{(m^k)}={k\choose
2}\left(m+\frac{1}{2}+\frac{1}{2(2m-1)}-\frac{4m}{{2m\choose
m}}\right).
\end{align*}
  Note that, when $m=2$, we recover already known results (see e.g.~\cite{FlN,Ri}).

\vspace{0.5cm}
\textbf{Acknowledgments.} The author would like to thank Christian
Krattenthaler and Einar Steingr\'\i msson for their encouragement 
during the preparation of this paper, and the anonymous referee for his helpful suggestions 
concerning the presentation of this paper.
Special acknowledgement is given to Christian Krattenthaler for a careful proofreading of this paper. 
All remaining errors are the author's.

%
%


\end{document}